\documentclass[12pt]{amsart}
\usepackage[vmargin=2.6cm,hmargin=2.2cm]{geometry} 
\usepackage[colorlinks,citecolor=blue,urlcolor=blue,bookmarks=false,hypertexnames=true]{hyperref}

\usepackage{enumerate}
\usepackage{color}
\usepackage{amsthm}
\usepackage{graphicx}
\usepackage{times}
\usepackage{xcolor}
\graphicspath{ {./images/} }
\usepackage{tikz}
\usepackage{multirow}
\usepackage{caption}
\usepackage{subcaption}
\usepackage[square,numbers]{natbib}
\usepackage{comment}
\usepackage[title]{appendix}
\usepackage{fancyhdr}	

\allowdisplaybreaks[3] 

\numberwithin{equation}{section}
\newtheorem{theorem}{Theorem}[section]

\newtheorem{lemma}{Lemma}[section]

\newtheorem{remark}{Remark}[section]

\newtheorem{assumption}{Assumption}[section]
\newtheorem{definition}{Definition}[section]

\newcommand{\dd}{\operatorname{d}\! }
\newcommand{\dt}{\operatorname{d}\! t}
\newcommand{\de}{\operatorname{d}\! e}
\newcommand{\ds}{\operatorname{d}\! s}
\newcommand{\dr}{\operatorname{d}\! r}

\newcommand{\dz}{\operatorname{d}\! z}

\newcommand{\dw}{\operatorname{d}\! W}

\newcommand{\cM}{\ensuremath{\mathcal{M}}}
\newcommand{\cZ}{\ensuremath{\mathcal{Z}}}
\newcommand{\cE}{\ensuremath{\mathcal{E}}}
\newcommand{\argmin}{\ensuremath{\operatorname*{argmin}}}

\newcommand{\esssup}{\ensuremath{\operatorname*{ess\;sup}}}

\newcommand{\nn}{\nonumber}

\newcommand{\E}{\mathbb{E}}
\newcommand{\R}{\mathbb{R}}
\newcommand{\pf}{\noindent\textbf{Proof:} }
\newcommand{\eof}{\hfill{$\Box$}}

\newcommand{\dpt}{\dd\mathbb{P}\otimes \dt\textrm{-a.e.}}
\newcommand{\dptv}{\dd\mathbb{P}\otimes \dt\otimes\dd\nu\textrm{-a.e.}}

\title[LQ with controlled jump size]{Constrained stochastic linear quadratic control under regime switching with controlled jump size}
\author[Shi]{Xiaomin Shi}
\author[Xu]{Zuo Quan Xu}
\date{\today}

\address{X.~Shi: School of Statistics and Mathematics, Shandong University of Finance and Economics, Jinan
250100, China.}
\email{{shixm@mail.sdu.edu.cn}}
\address{Z.Q.~Xu: Department of Applied Mathematics, The Hong Kong Polytechnic University, Kowloon, Hong Kong, China.}
\email{{maxu@polyu.edu.hk}}

\begin{document}
\setcitestyle{numbers}
\maketitle


\begin{abstract}
In this paper, we examine a stochastic linear-quadratic control problem characterized by regime switching and Poisson jumps. All the coefficients in the problem are random processes adapted to the filtration generated by Brownian motion and the Poisson random measure for each given regime. The model incorporates two distinct types of controls: the first is a conventional control that appears in the continuous diffusion component, while the second is an unconventional control, dependent on the variable $z$, which influences the jump size in the jump diffusion component. Both controls are constrained within general closed cones. By employing the Meyer-It\^o formula in conjunction with a generalized squares completion technique, we rigorously and explicitly derive the optimal value and optimal feedback control. These depend on solutions to certain multi-dimensional fully coupled stochastic Riccati equations, which are essentially backward stochastic differential equations with jumps (BSDEJs). We establish the existence of a unique nonnegative solution to the BSDEJs. One of the major tools used in the proof is the newly established comparison theorems for multidimensional BSDEJs.

\end{abstract}

\subsection*{Keywords:}linear-quadratic control; regime switching; controlled jump size; fully coupled stochastic Riccati equations; backward stochastic differential equations with jumps.
\subsection*{Mathematics Subject Classification (2020):} 60H30. 60J28. 60J76. 93E20.

\section{Introduction}

Since the pioneering work of Wonham \cite{Wo}, stochastic linear-quadratic (LQ) theory has been extensively studied by numerous researchers. For instance, Bismut \cite{Bi} was the first one who studied stochastic LQ problems with random coefficients. In order to obtain the optimal random feedback control, he formally derived a stochastic Riccati equation (SRE). But he could not solve the SRE in the general case. It is Kohlmann and Tang \cite{KT}, for the first time, that established the existence and uniqueness of the one-dimensional SRE.
Tang \cite{Tang03, Tang15} made another breakthrough and proved the existence and uniqueness of the matrix valued SRE with uniformly positive control weighting matrix using two different approaches. Sun, Xiong and Yong \cite{SXY} studied the indefinite stochastic LQ problem with random coefficients. Hu and Zhou \cite{HZ} solved the stochastic LQ problem with cone control constraint.
Zhang, Dong and Meng \cite{ZDM} made a great progress in solving stochastic LQ control and related SRE with jumps with uniformly definite control weight by inverse flow technique. Li, Wu and Yu \cite{LWY} considered the stochastic LQ problem with jumps in the indefinite case. Please refer to Chapter 6 in Yong and Zhou \cite{YZ} for a systematic account on this subject.

Stochastic LQ problems for Markovian regime switching system were studied in
Wen, Li and Xiong \cite{WLX} and Zhang, Li and Xiong \cite{ZLX} where weak closed-loop solvability, open-loop solvability and closed-loop solvability were established. But the coefficients are assumed to be \emph{deterministic} functions of time $t$ for each given regime $i$ in the above papers,
so their SREs are indeed deterministic ordinary differential equations (ODEs).
Hu, Shi and Xu \cite{HSX, HSX2} formulated cone-constrained stochastic LQ problems with regime switching on finite time horizon and infinite time horizon respectively, in which the coefficients are \emph{stochastic} processes adapted to the filtration generated by the Brownian motion for each give regime $i$.
Due to the randomness of the coefficients, the corresponding SREs in \cite{HSX, HSX2} are actually BSDEs.

In this paper, we generalize the LQ problem in \cite{HSX} to a model in which the coefficients are \emph{stochastic} processes adapted to the filtration generated by the Brownian motion and the Poisson random measure for each give regime $i$. In addition to a usual control $u_1$, we introduce a second control $u_2(z)$ depending on the jump size $z$. The motivations to incorporating the second control are, in insurance area, the optimal reinsurance strategies may depend on the claim size in general, see, e.g., Liu and Ma \cite{LM} and Wu, Shen, Zhang and Ding \cite{WSZD}; and in controllability issues for stochastic systems with jump diffusions, a control depending on the jump size is necessary as a consequence of martingale representation theorem of Poisson random measures, see, e.g., Goreac \cite{Goreac} and Song \cite{Song}. The application of this kind of stochastic LQ model in a optimal liquidation problem with dark pools can be found in our working paper \cite{FSX}.

The first main contribution of this paper is to provide a pure analysis method (using tools like approximation technique, comparison theorem for multi-dimensional BSDEJs, log transformation, etc) of the existence of a unique solution to the corresponding system of SREs, which is a $2\ell$-dimensional coupled BSDEJs. This is interesting in its own right from the point of view of BSDE theory.
Note that even though the SREs in \cite{HSX} are $2\ell$-dimensional, they are partially coupled, that is, the first $\ell$ equations for $\{P^i_1\}_{i\in\cM}$ and the second $\ell$ equations for $\{P^i_2\}_{i\in\cM}$ are totally decoupled. But in our new model, the equation for $P^i_1$ also depends on $P^i_2,\Gamma^i_2$, rendering the $2\ell$-dimensional SREs in our new model are fully coupled.
This more complicated phenomenon comes from the fact that, to the best of our knowledge, the optimal state process will probably change its sign at the jump time of the underlying Poisson random measure. Compared with the $2$-dimensional SREs in Hu, Shi and Xu \cite{HSX4}, here we need to study $2\ell$-dimensional SREs because of the new coupling terms $\sum_j q^{ij}P^j_1$ and $\sum_j q^{ij}P^j_2$.
The second main contribution is to give a rigorous verification theorem of the optimal value and optimal control, using the unique solution to the corresponding system of SREs, Meyer-It\^o's formula, a generalized squares completion technique and some delicate analysis.

The rest part of this paper is organized as follows. In Section \ref{section:fm}, we formulate a constrained stochastic LQ control problem with regime switching, controlled jump size and random coefficients. Section \ref{section:Ri} is devoted to proving the existence of a unique nonnegative solution to the related $2\ell$-dimensional fully coupled SREs in standard and singular cases.
In Section \ref{section:Veri}, we solve the LQ problem by establishing a rigorous verification theorem.

\section{Problem formulation}\label{section:fm}
Let $(\Omega, \mathcal F, \mathbb{F}, \mathbb{P})$ be a fixed complete filtered probability space. The filtration $\mathbb{F}=\{\mathcal F_t, t\geq0\}$ is generated by the following three independent random sources augmented by all the $\mathbb{P}$-null sets. \medskip
\begin{itemize}
\item The first random source is a standard $n$-dimensional Brownian motion $W_t=(W_{1,t}, \ldots, W_{n_1,t})^{\top}$.\medskip

\item The second one is an $n_2$-dimensional Poisson random measure $N=(N_1, \ldots,N_{n_2})^{\top}$ defined on $\R_+\times\cZ$, where $\mathcal{Z}\subset\R^{\ell}\setminus\{0\}$ is a nonempty Borel subset of some Euclidean space. For each $k=1,\ldots,n_2$, $N_k$ posses the same stationary compensator (intensity measure) $\nu(\dz)\dt$ satisfying $\nu(\cZ)<\infty$. The compensated Poisson random measure is denoted by $\tilde N(\dt,\dz)$.\medskip

\item The third one is a continuous-time stationary Markov chain $\alpha_t$ valued in a finite state space $\mathcal M=\{1, 2, \ldots, \ell\}$ with $\ell\geq 1$. The Markov chain has a generator $Q=(q_{ij})_{\ell\times\ell}$ with $q_{ij}\geq0$ for $i\neq j$ and $\sum_{j=1}^{\ell}=0$ for every $i\in\cM$.
\end{itemize}
Besides the filtration $\mathbb{F}$, we will often use the filtration $\mathbb{F}^{W,N}=\{\mathcal F^{W,N}_t, t\geq0\}$ which is generated by the Brownian motion $W$ and the Poisson random measures $N$ and augmented by all the $\mathbb{P}$-null sets. Throughout the paper, let $T$ denote a fixed positive constant, $\mathcal{P}$ (resp. $\mathcal{P}^{W,N}$) denote the $\mathbb{F}$ (resp. $\mathbb{F}^{W,N}$)-predictable $\sigma$-field on $\Omega\times[0,T]$, and $\mathcal{B}(\cZ)$ denote the Borelian $\sigma$-field on $\cZ$.

We denote by $\R^\ell$ the set of $\ell$-dimensional column vectors, by $\R^\ell_+$ the set of vectors in $\R^\ell$ whose components are nonnegative, by $\R^{\ell\times n}$ the set of $\ell\times n$ real matrices, by $\mathbb{S}^n$ the set of $n\times n$ symmetric real matrices, by $\mathbb{S}^n_+$ the set of $n\times n$ nonnegative definite real matrices, and by $\mathbf{1}_n$ the $n$-dimensional identity matrix. For any vector $Y$, we denote $Y_i$ as its $i$-th component.
For any matrix $M=(m_{ij})$, we denote its transpose by $M^{\top}$, and its norm by $|M|=\sqrt{\sum_{ij}m_{ij}^2}$. If $M\in\mathbb{S}^n$ is positive definite (resp. positive semidefinite), we write $M>$ (resp. $\geq$) $0.$ We write $A>$ (resp. $\geq$) $B$ if $A, B\in\mathbb{S}^n$ and $A-B>$ (resp. $\geq$) $0.$
We write the positive and negative parts of $x\in\R$ as $x^+=\max\{x, 0\}$ and $x^-=\max\{-x, 0\}$ respectively. The elementary inequality $|a^{\top}b|\leq c|a|^2+\frac{|b|^2}{2c}$ for any $a,b\in\R^{n}$, $c>0$, will be used frequently without claim.
Throughout the paper, we use $c$ to denote a suitable positive constant, which is independent of $(t,\omega, i)$ and can be different from line to line.

\subsection*{Notation}
We use the following notation throughout the paper:
\begin{align*}
L^{\infty}_{\mathcal{F}_T}(\Omega;\R)&=\Big\{\xi:\Omega\rightarrow
\R\;\Big|\;\xi\mbox { is }\mathcal{F}_{T}\mbox{-measurable, and essentially bounded}\Big\}, \\
L^{2}_{\mathbb{F}}(0, T;\R)&=\Big\{\phi:[0, T]\times\Omega\rightarrow
\R\;\Big|\;\phi\mbox{ is } \mathbb{F}%
\mbox{-predictable and }\E\int_{0}^{T}|\phi_t|^{2}\dt<\infty
\Big\}, \\
L^{\infty}_{\mathbb{F}}(0, T;\R)&=\Big\{\phi:[0, T]\times\Omega
\rightarrow\R\;\Big|\;\phi\mbox{ is }\mathbb{F}%
\mbox{-predictable and essentially bounded} \Big\},\\
L^{2,\nu}(\R)&=\Big\{\phi:\cZ\rightarrow\R\mbox{ is measurable with} \ \|\phi(\cdot)\|^2_{\nu}:=\int_{\cZ}\phi(z)^2\nu(\dz)<\infty\Big\},\\
L^{\infty,\nu}(\R)&=\Big\{\phi:\cZ\rightarrow\R\mbox{ is measurable and} \ \phi \mbox{ is bounded }\dd\nu \textrm{-a.e.}\Big\},\\
L^{2,\nu}_{\mathcal{P}}(0, T;\R)&=\Big\{\phi:[0, T]\times\Omega\times\cZ\rightarrow
\R\;\Big|\;\phi\mbox{ is } \mathcal{P}\otimes\mathcal{B}(\cZ)%
\mbox{-measurable }\\
&\qquad\mbox{ \ \ \ \ and }\E\int_{0}^{T}\int_{\cZ}|\phi_t(z)|^{2}\nu(\dz)\dt<\infty
\Big\},\\
L^{\infty}_{\mathcal{P}}(0, T;\R)&=\Big\{\phi:[0, T]\times\Omega\times\cZ
\rightarrow\R\;\Big|\;\phi \mbox{ is }\mathcal{P}\otimes\mathcal{B}(\cZ)%
\mbox{-measurable and essentially bounded} \Big\},\\
S^{\infty}_{\mathbb{F}}(0,T;\R)&=\Big\{\phi:\Omega\times[0,T]\to \R
\;\Big|\;\phi \mbox{ is c\`ad-l\`ag, $\mathbb{F}$-adapted and essentially bounded}\Big\}.
\end{align*}
These definitions are generalized in the obvious way to the cases that $\mathcal{F}$ is replaced by $\mathcal F^{W,N}$, $\mathbb{F}$ by $\mathbb{F}^{W,N}$, $\mathcal{P}$ by $\mathcal{P}^{W,N}$ and $\R$ by $\R^n$, $\R^{n\times m}$ or $\mathbb{S}^n$.
In our argument, $t$, $\omega$, ``almost surely'' and ``almost everywhere'', will be suppressed for simplicity in many circumstances, when no confusion occurs.
All the processes and maps considered in this paper, unless otherwise stated, are stochastic, so, for notation simplicity, we will not write their dependence on $\omega$ explicitly.
Equations and inequalities shall be understood to hold true $\dptv$
For a random variable or stochastic process $X$, we write $X\gg1$ (resp. $X\ll1$) if there exists a constant $c>0$ such that $X\geq c$ (resp. $|X|\leq c$).

Consider the following real-valued linear stochastic differential equation (SDE) with jumps:
\begin{align}
\label{state}
\begin{cases}
\dd X_t=\left[A_t^{\alpha_{t-}}X_{t-}+(B^{\alpha_{t-}}_{1,t})^{\top}u_{1,t}+\int_{\cZ}B^{\alpha_{t-}}_{2,t}(z)^{\top}u_{2,t}(z)\nu(\dz)\right]\dt\\
\qquad\qquad\qquad+\left[C^{\alpha_{t-}}_tX_{t-}+D^{\alpha_{t-}}_t u_{1,t}\right]^{\top}\dw_t\\
\qquad\qquad\qquad+\int_{\cZ}\left[E_t^{\alpha_{t-}}(z)X_{t-}+F_t^{\alpha_{t-}}(z)u_{2,t}(z)\right]^{\top}\tilde N(\dt,\dz), \quad t\in[0,T], \\
X_0=x,~~ \alpha_0=i_0,
\end{cases}
\end{align}
where $A^{i}, \ B_1^{i}, \ C^{i}, \ D^{i}$ are all $\mathbb{F}^{W,N}$-predictable processes, and $B^{i}_{2}(\cdot), \ E^{i}(\cdot), \ F^{i}(\cdot)$ are $\mathcal{P}^{W,N}\otimes\mathcal{B}(\cZ)$-measure processes of suitable sizes, $(u_1,u_2)$ is the control and $x\in\R$, $i_0\in\cM$ are the known initial values.

Let $\Pi_1$, $\Pi_2$ be two given closed cones (not necessarily convex) in $\R^{m_1}$ and $\R^{m_2}$, respectively.
The class of admissible controls is defined as the set
\begin{align*}
\mathcal{U} &:=\Big\{(u_{1},u_{2}) \;\Big|\;u_{1}\in L^2_\mathbb{F}(0, T;\R^{m_1}),~u_{1,t}\in\Pi_1, ~\dpt, \\
&\qquad\qquad\qquad\ \mbox{and }\ u_{2}\in L^{2,\nu}_{\mathcal{P}}(0, T;\R^{m_2}),~ u_{2,t}\in\Pi_2, ~\dptv
\Big\}.%
\end{align*}
If $u\equiv (u_{1},u_{2})\in\mathcal{U}$, then the SDE \eqref{state} admits a unique strong solution $X$, and we refer to $(X, u)$ as an admissible pair.

Let us now state our stochastic linear quadratic optimal control problem as follows:
\begin{align}
\begin{cases}
\mbox{Minimize} &\ J(u;x,i_0 )\smallskip\\
\mbox{subject to} &\ u \in\mathcal{U},
\end{cases}
\label{LQ}%
\end{align}
where the cost functional $J$ is given as the following quadratic form
\begin{align}\label{costfunc}
J(u;x,i_0) &:=\E\Big\{ G^{\alpha_{T}}X_T^2
+\int_0^T\Big[u_{1,t}^{\top}R^{\alpha_{t}}_{1,t}u_{1,t}
+Q^{\alpha_{t}}_tX_t^2+\int_{\cZ}u_{2,t}(z)^{\top}R^{\alpha_{t}}_{2,t}(z)u_{2,t}(z)\nu(\dz)\Big]\dt\Big\}.
\end{align}
The optimal value of the problem is defined as
\begin{align*}
V(x,i_0)=\inf_{u \in\mathcal{U}}~J( u;x,i_0).
\end{align*}
Problem \eqref{LQ} is said to be solvable, if there exists a control $u^*\in\mathcal{U}$ such that
\begin{align*}
-\infty<J( u^*;x,i_0)\leq J( u;x,i_0), \quad \forall\; u\in\mathcal{U},
\end{align*}
in which case, $u^*$ is called an optimal control for problem \eqref{LQ} and one has
\begin{align*}
V(x,i_0)= J( u^*;x,i_0).
\end{align*}

\begin{remark}
By choosing $\Pi_1=\{0\}$ (resp. $\Pi_2=\{0\}$), our model covers the case of $R_1^i=0$, $D^i=0$, $B^i_1=0$ (resp. $R_2^i=0$, $F^i=0$, $B^i_2=0$). In particular, our model covers the pure jump (i.e. $(B^i_1,C^i,D^i,R_1^i)=0$) and pure diffusion (i.e., $(B^i_2,E^i,F^i,R_2^i)=0$) models.
\end{remark}

Throughout this paper, we put the following assumption on the coefficients.
\begin{assumption} \label{assu1}
It holds, for every $i\in\cM$, that
\begin{align*}
\begin{cases}
A^{i}\in L_{\mathbb{F}^{W,N}}^\infty(0, T;\R), \
B^{i}_1\in L_{\mathbb{F}^{W,N}}^\infty(0, T;\R^{m_1}), \
B^{i}_2\in L_{\mathcal{P}^{W,N}}^\infty(0,T;\R^{m_2}),\\
C^{i} \in L_{\mathbb{F}^{W,N}}^\infty(0, T;\R^{n_1}), \
D^{i}\in L_{\mathbb{F}^{W,N}}^\infty(0, T;\R^{n_1\times m_1}), \\
E^{i}\in L_{\mathcal{P}^{W,N}}^\infty(0,T;\R^{n_2}),\
F^{i}\in L_{\mathcal{P}^{W,N}}^\infty(0,T;\R^{n_2\times m_2}),\\
R^{i}_1\in L_{\mathbb{F}^{W,N}}^\infty(0, T;\mathbb{S}^{m_1}_+), \
R^{i}_2 \in L_{\mathcal{P}^{W,N}}^\infty(0,T;\mathbb{S}^{m_2}_+),\\
Q^{i}\in L_{\mathbb{F}^{W,N}}^\infty(0, T;\R_+), \
G^{i}\in L_{\mathcal{F}^{W,N}_T}^\infty(\Omega;\R_+).
\end{cases}
\end{align*}
\end{assumption}
%
%
%
Under Assumption \ref{assu1}, the cost functional \eqref{costfunc} is nonnegative, hence problem \eqref{LQ} is well-posed.
\bigskip
%

Besides Assumption \ref{assu1}, we need the following hypothesizes:
\begin{enumerate}[1.]
\item \label{R11} $R^{i}_1\geq \delta \mathbf{1}_{m_1}$.
\item \label{R12} $(D^{i})^{\top}D^i\geq \delta \mathbf{1}_{m_1}$.
\item \label{R21} $R^{i}_2\geq \delta \mathbf{1}_{m_2}$.
\item \label{R22} $(F^{i})^{\top}F^i\geq \delta \mathbf{1}_{m_2}$.
\end{enumerate}
We will consider the problem under one of following two assumptions.
\begin{assumption} [Standard case]
\label{assu2} There exists a constant $\delta>0$ such that
both hypothesizes \ref{R11} and \ref{R21} hold.
\end{assumption}

\begin{assumption}[Singular case]
\label{assu3}
There exists a constant $\delta>0$ such that $G^i\geq\delta$ and one of the following hold:
\begin{enumerate}
\item[Case I.] Both hypothesizes \ref{R12} and \ref{R21} hold;
\item[Case II.] Both hypothesizes \ref{R12} and \ref{R22} hold;
\item[Case III.] Both hypothesizes \ref{R11} and \ref{R22} hold.
\end{enumerate}
\end{assumption}

\section{Solvability of the Riccati equations}\label{section:Ri}
For any $i\in\cM$, $j=1,2,$ let us denote by $E^i_k$ the $k$-th component of $E^i$, $\Gamma^i_{jk}$ the $k$-th component of $\Gamma^i_j$ and $F^i_k$ the $k$-th row of $F^i$, $k=1,...,n_2$.
To solve problem \eqref{LQ}, we need to study the following $2\ell$-dimensional SRE with jumps:
\begin{align}
\label{P}
\begin{cases}
\dd P_{1,t}^i=-\Big[(2A^i+|C^i|^2)P^i_{1,t-}+2(C^i)^{\top}\Lambda^i_1+Q^i+H_{11}^{i,*}(P^i_1,\Lambda^i_1)\\
\qquad\qquad\qquad\qquad+\int_{\cZ}H_{12}^{i,*}(z,P^i_1,P^i_2,\Gamma^i_1, \Gamma^i_2)
\nu(\dz)+\sum_{j=1}^{\ell}q^{ij}P_{1}^j\Big]\dt\\
\qquad\quad\;+(\Lambda^i_1)^{\top}\dw+\int_{\cZ}\Gamma^i_1(z)^{\top}\tilde N(\dt,\dz),\bigskip\\
\dd P_{2,t}^i=-\Big[(2A^i+|C^i|^2)P^i_{2,t-}+2(C^i)^{\top}\Lambda^i_2+Q^i+H_{21}^{i,*}(P^i_2,\Lambda^i_2)\\
\qquad\qquad\qquad\qquad+\int_{\cZ}H_{22}^{i,*}(z,P^i_1,P^i_2,\Gamma^i_1, \Gamma^i_2)
\nu(\dz)+\sum_{j=1}^{\ell}q^{ij}P_{2}^j\Big]\dt\\
\qquad\quad\;+(\Lambda^i_2)^{\top}\dw+\int_{\cZ}\Gamma^i_2(z)^{\top}\tilde N(\dt,\dz),\bigskip\\
P^i_{1,T}=G, \ P^i_{2,T}=G, ~~
R^i_{1,t}+P^i_{1,t}(D^i_t)^{\top}D^i_t> 0, \ R^i_{1,t}+P^i_{2,t}(D^i_t)^{\top}D^i_t> 0, ~~i\in\cM,\\
\end{cases}
\end{align}
where, for any $(z,P_1,P_2,\Lambda,\Gamma_1,\Gamma_2)\in\cZ\times\R_+\times\R_+\times\R^n\times \R \times \R$,
\begin{align*}
H^i_{11}(v,P_1,\Lambda)&:=v^{\top}(R^i_1+P_1(D^i)^{\top}D^i)v+2(P_1(B^i_1+(D^i)^{\top}C^i)+(D^i)^{\top}\Lambda)^{\top}v \\
H^i_{12}(v,z,P_1,P_2,\Gamma_1, \Gamma_2)&:=v^{\top}R^i_2v+\sum_{k=1}^{n_2}(P_1+\Gamma_{1k})\big[((1+E^i_k+F^i_k v)^+)^2-1\big]\\
&\quad-2P_1\sum_{k=1}^{n_2}(E^i_k+F^i_k v)+2P_1(B^i_2)^{\top}v+\sum_{k=1}^{n_2}(P_2+\Gamma_{2k}) ((1+E^i_k+F^i_k v)^-)^2, \\
H^i_{21}(v,P_2,\Lambda)&:=v^{\top}(R^i_1+P_2(D^i)^{\top}D^i)v-2(P_2(B^i+(D^i)^{\top}C^i)+(D^i)^{\top}\Lambda)^{\top}v,\\
H^i_{22}(v,z,P_1,P_2,\Gamma_1, \Gamma_2)&:=v^{\top}R^i_2v+\sum_{k=1}^{n_2}(P_2+\Gamma_{2k})\big[((-1-E^i_k+F^i_k v)^-)^2-1\big]\\
&\quad-2P_2\sum_{k=1}^{n_2}(E^i_k-F^i_kv)-2P_2(B^i_2)^{\top}v+\sum_{k=1}^{n_2}(P_1+\Gamma_{1k})((-1-E^i_k+F^i_k v)^+)^2,
\end{align*}
and
\begin{align*}
H_{11}^{i,*}(t,P_1,\Lambda)&:=\inf_{v\in\Pi_1}
H_{11}^i(v,P_1,\Lambda),\\
H_{12}^{i,*}(t,z,P_1,P_2,\Gamma_1, \Gamma_2)&:=\inf_{v\in\Pi_2}
H_{12}^i(v,z,P_1,P_2,\Gamma_1, \Gamma_2),\\
H_{21}^{i,*}(t,P_2,\Lambda)&:=\inf_{v\in\Pi_1}
H_{21}^i(v,P_2,\Lambda),\\
H_{22}^{i,*}(t,z,P_1,P_2,\Gamma_1, \Gamma_2)&:=\inf_{v\in\Pi_2}
H_{22}^i(v,z,P_1,P_2,\Gamma_1, \Gamma_2).
\end{align*}
Noting in above, we omit the argument $t$.
Because the generators in \eqref{P} depend on all $P^i_j$s, hence \eqref{P} is a fully coupled BSDE with jumps.



\begin{definition}\label{def}
A vector of stochastic process $(P_j^i,\Lambda_j^i,\Gamma_j^i)_{i\in\cM,\;j=1,2}$ is called a solution to the BSDEJ \eqref{P} if it satisfies all the equations and constraints in \eqref{P}, and $(P_j^i,\Lambda_j^i,\Gamma_j^i)\in S^{\infty}_{\mathbb{F}^{W,N}}(0,T;\R)\times L^{2}_{\mathbb{F}^{W,N}}(0,T;\R^{n_1})\times L^{\infty,\nu}_{\mathcal{P}^{W,N}}(0, T;\R^{n_2})$ for all $i\in\cM$, $j=1,2$. Furthermore, the solution is called nonnegative if $P^i_j\geq0$, $P^i_j+\Gamma^i_j\geq0$, and called uniformly positive if $P^i_j\gg1$ and $P^i_j+\Gamma^i_j\gg1$, for all $i\in\cM$, $j=1,2$.
\end{definition}

Before giving the proof of main theorem in this section, let us recall the definition of bounded mean oscillation martingales, briefly called BMO martingales. Please refer to Kazamaki \cite{Ka} for a systematic account on continuous BMO martingales. A process $\int_0^{t}\phi_s^\top dW_s$ is called a BMO martingale if and only if there exists a constant $c>0$ such that
\[
\E\Big[\int_{\tau}^T|\phi_s|^2\ds\;\Big|\;\mathcal{F}^{W,N}_{\tau}\Big]\leq c
\]
for all $\mathbb{F}^{W,N}$ stopping times $\tau\leq T$.
In the uniqueness part of the Theorem \ref{existence}, we will use the following property of BMO martingales:
If $\int_0^{t}\phi_s^\top dW_s$ is a BMO martingale on $[0,T]$, then the Dol\'eans-Dade stochastic exponential $\mathcal{E}\big(\int_0^{t}\phi_s^{\top}dW_s\big)$ is a uniformly integrable martingale on $[0,T]$.

The following comparison theorem for multi-dimensional BSDEJs was firstly established in \cite[Theorem 2.2]{HSX4}. We list it here as it plays crucial role in the solvability of the BSDEJ \eqref{P}.
\begin{lemma}
\label{comparison}
Suppose, for every $ i\in \{1,2,...,m\}$,
\begin{align*}
(Y_i, Z_i, \Phi_i), (\overline Y_i, \overline Z_i, \overline\Phi_i)\in S^{2}_{\mathbb{F}}(0,T;\R)\times L^{2}_{\mathbb F}(0, T;\R^{n})\times L^{2,\nu}_{\mathcal{P}}(0, T;\R),
\end{align*}
and they satisfy BSDEJs
\begin{align*}
Y_{i,t}&=\xi_i+\int_t^T f_i(s,Y_{s-}, Z_{i,s}, \Phi_{s})\ds-\int_t^T Z_{i,s}^{\top}\dw_s-\int_t^T\int_{\cE}\Phi_{i,s}(e) \widetilde N(\ds,\de),
\end{align*}
and
\begin{align*}
\overline Y_{i,t}&=\overline\xi_i+\int_t^T \overline f_i(s, \overline Y_{s-}, \overline Z_{i,s}, \overline \Phi_{s})\ds-\int_t^T\overline Z_{i,s}^{\top}\dw_s-\int_t^T\int_{\cE}\overline\Phi_{i,s}(e) \widetilde N(\ds,\de).
\end{align*}
Also suppose that, for all $i\in\{1,2,...,m\}$ and $s\in[0,T]$,
\begin{enumerate}[(1)]
\item
\label{cond-boundary}
$\xi_i\leq\overline\xi_i$;
\item \label{cond-gamma}
there exists a constant $c>0$ such that
\begin{align*}
&\quad\;f_i(s,Y_{s-}, Z_{i,s}, \Phi_{1,s}, \cdots, \Phi_{i,s}, \cdots, \Phi_{\ell,s})\\
&\qquad\quad-f_i(s,Y_{s-}, Z_{i,s}, \Phi_{1,s}, \cdots, \overline\Phi_{i,s}, \cdots, \Phi_{\ell,s}) \\
&\leq c \int_{\cE} (\Phi_{i,s}(e)-\overline \Phi_{i,s}(e))^{+}\nu(\de)
+\int_{\cE} |\Phi_{i,s}(e)-\overline \Phi_{i,s}(e)|\nu(\de);
\end{align*}
\item \label{cond-growth} there exists a constant $c>0$ such that
\begin{align*}
&\quad\;f_i(s,Y_{s-}, Z_{i,s}, \Phi_{1,s}, \cdots, \overline\Phi_{i,s}, \cdots, \Phi_{\ell,s})-f_i(s, \overline Y_{s-}, \overline Z_{i,s}, \overline \Phi_{s})\\
&\leq c\Big (| Y_{i,s-}- \overline Y_{i,s-}|+\sum_{j\neq i} (Y_{j,s-}- \overline Y_{j,s-})^{+}+|Z_{i,s}- \overline Z_{i,s}|\\
&\qquad\quad+\sum_{j\neq i} \int_{\cE}( Y_{j,s-}+\Phi_{j,s}(e)-\overline Y_{j,s-}
-\overline\Phi_{j,s}(e))^+\nu(\de)\Big);
\end{align*}
\item $f_{i}(\cdot, 0,0,0)$ and $\overline f_{i}(\cdot, 0,0,0)\in L^{2}_{\mathbb F}(0, T;\R)$;
\item both $f_{i}$ and $\overline f_i$ are Lipschitz in $(y,z,\phi)$;~\text{and}
\item \label{cond-size} $f_i(s, \overline Y_{s-}, \overline Z_{i,s}, \overline \Phi_{s})
\leq \overline f_i(s, \overline Y_{s-}, \overline Z_{i,s}, \overline \Phi_{s})$.
\end{enumerate}
Then $\mathbb{P}\{Y_{i,t}\leq\overline Y_{i,t}, \forall t\in[0,T]\}=1$ for all $i\in\{1,2,...,m\}$.
\end{lemma}

\begin{theorem}\label{existence}
Under Assumptions \ref{assu1} and \ref{assu2}, the BSDEJ \eqref{P} admits a unique nonnegative
solution $(P_j^i,\Lambda_j^i,\Gamma_j^i)_{i\in\cM,\;j=1,2}$.
\end{theorem}
\pf
For each natural number $k$, define maps
\begin{align*}
H_{11}^{i,*,k}(P_{1},\Lambda) &:=\inf_{v\in\Pi_1,|v|\leq k}
H_{11}^{i}(v,P_{1},\Lambda),\\
H_{12}^{i,*,k}(z,P_{1},P_{2},\Gamma_{1},\Gamma_2) &:=\inf_{v\in\Pi_2,|v|\leq k}
H_{12}^{i}(v,z,P_{1},P_{2},\Gamma_{1},\Gamma_2),\\
H_{21}^{i,*,k}(P_{2},\Lambda) &:=\inf_{v\in\Pi_1,|v|\leq k}
H_{21}^{i}(v,P_{2},\Lambda),\\
H_{22}^{i,*,k}(z,P_{1},P_{2},\Gamma_{1},\Gamma_2) &:=\inf_{v\in\Pi_2,|v|\leq k}
H_{22}^{i}(v,z,P_{1},P_{2},\Gamma_{1},\Gamma_2).
\end{align*}
Then they are uniformly Lipschitz in $(P_{1},P_{2},\Lambda,\Gamma_1,\Gamma_2)$ and decreasingly approach to
$H_{11}^{i,*}$, $H_{12}^{i,*}$, $H_{21}^{i,*}$, $H_{22}^{i,*}$ and respectively as $k$ goes to infinity.

For each $k$,
the following BSDE
\begin{align}
\label{Ptrun}
\begin{cases}
\dd P_{1,t}^{i,k}=-\Big[(2A^i+|C^i|^2)P^{i,k}_{1,t-}+2(C^i)^{\top}\Lambda^{i,k}_1+Q^i+
H_{11}^{i,*,k}(P^{i,k}_1,\Lambda^{i,k}_1)+\sum_{j=1}^{\ell}q^{ij}P_{1}^{j,k}\\
\qquad\qquad\qquad\qquad+\int_{\cZ}H_{12}^{i,*,k}(P^{i,k}_1,P^{i,k}_2,\Gamma^{i,k}_1, \Gamma^{i,k}_2)
\nu(\dz)\Big]\dt\\
\qquad\quad\;+(\Lambda^{i,k}_1)^{\top}\dw+\int_{\cZ}\Gamma^{i,k}_1(z)^{\top}\tilde N(\dt,\dz),\bigskip\\
\dd P_{2,t}^{i,k}=-\Big[(2A^i+|C^i|^2)P^{i,k}_{2,t-}+2(C^i)^{\top}\Lambda^{i,k}_2+Q^i+
H_{21}^{i,*,k}(P^{i,k}_2,\Lambda^{i,k}_2)+\sum_{j=1}^{\ell}q^{ij}P_{2}^{j,k}\\
\qquad\qquad\qquad\qquad+\int_{\cZ}H_{22}^{i,*,k}(P^{i,k}_1,P^{i,k}_2,\Gamma^{i,k}_1, \Gamma^{i,k}_2)\nu(\dz)\Big]\dt\\
\qquad\quad\;+(\Lambda^{i,k}_2)^{\top}\dw+\int_{\cZ}\Gamma^{i,k}_2(z)^{\top}\tilde N(\dt,\dz),\bigskip\\
P^{i,k}_{1,T}=G, \ P^{i,k}_{2,T}=G, \ i\in\cM,
\end{cases}
\end{align}
is a $2\ell$-dimensional BSDEJ with a Lipschitz generator. According to \cite[Lemma 2.4]{TL}, it admits a unique solution $(P_j^{i,k},\Lambda_j^{i,k},\Gamma_j^{i,k})_{i\in\cM,\;j=1,2}$
such that $$(P_j^{i,k},\Lambda_j^{i,k},\Gamma_j^{i,k})\in S^{2}_{\mathbb{F}^{W,N}}(0, T;\R)\times L^{2}_{\mathbb{F}^{W,N}}(0, T;\R^{n_1})\times L^{2,\nu}_{\mathcal{P}^{W,N}}(0, T;\R^{n_2}), \ i\in\cM, \ j=1,2.$$

Next we will show that $(P_1^{i,k},P_2^{i,k})_{i\in\cM}$ are lower and upper bounded.
Actually, the following two linear (with bounded coefficients) BSDEJs (see, e.g., \cite[Proposition 2.2]{BBP})
\begin{align}
\label{overlineP}
\begin{cases}
\dd \overline P_{1,t}^i=-\Big[(2A^i+|C^i|^2)\overline P^i_{1,t-}+2(C^i)^{\top}\overline \Lambda^i_1+Q^i+\int_{\cZ}H_{12}^{i}(0,\overline P^i_1,\overline P^i_2,\overline \Gamma^i_1, \overline \Gamma^i_2)
\nu(\dz)\\
\qquad\qquad+\sum_{j=1}^{\ell}q^{ij}\overline P_{1}^j\Big]\dt+(\overline \Lambda^i_1)^{\top}\dw+\int_{\cZ}\overline \Gamma^i_1(z)^{\top}\tilde N(\dt,\dz),\\
\dd \overline P^i_{2,t}=-\Big[(2A^i+|C^i|^2)\overline P^i_{2,t-}+2(C^i)^{\top}\overline \Lambda^i_2+Q^i+\int_{\cZ}H_{22}^{i}(0,\overline P^i_1,\overline P^i_2,\overline \Gamma^i_1, \overline \Gamma^i_2)
\nu(\dz)\\
\qquad\qquad+\sum_{j=1}^{\ell}q^{ij}\overline P_{2}^j\Big]\dt+(\overline \Lambda^i_2)^{\top}\dw+\int_{\cZ}\overline \Gamma^i_2(z)^{\top}\tilde N(\dt,\dz),\\
\overline P^i_{1,T}=G, \ \overline P^i_{2,T}=G, \ i\in\cM.
\end{cases}
\end{align}
and
\begin{align}
\label{underlineP}
\begin{cases}
\dd\underline P_{1,t}^i=-\Big[(2A^i+|C^i|^2)\underline P^i_{1,t-}+2(C^i)^{\top}\underline \Lambda^i_1+\sum_{j=1}^{\ell}q^{ij}\underline P_{1}^j\Big]\dt+(\underline \Lambda^i_1)^{\top}\dw+\int_{\cZ}\underline \Gamma^i_1(z)^{\top}\tilde N(\dt,\dz),\\
\dd\underline P_{2,t}^i=-\Big[(2A^i+|C^i|^2)\underline P^i_{2,t-}+2(C^i)^{\top}\underline \Lambda^i_2+\sum_{j=1}^{\ell}q^{ij}\underline P_{2}^j\Big]\dt+(\underline \Lambda^i_2)^{\top}\dw+\int_{\cZ}\underline \Gamma^i_2(z)^{\top}\tilde N(\dt,\dz),\\
\underline P^i_{1,T}=0, \ \underline P^i_{2,T}=0, \ i\in\cM.
\end{cases}
\end{align}
admit unique uniformly bounded solutions
$(\overline P_j^i, \overline\Lambda_j^i, \overline \Gamma_j^i)_{i\in\cM,\;j=1,2}$
and
$(\underline P_j^i, \underline\Lambda_j^i, \underline \Gamma_j^i)_{i\in\cM,\;j=1,2}$ respectively. Clearly, $(\underline P_j^i, \underline\Lambda_j^i, \underline \Gamma_j^i)_{i\in\cM,\;j=1,2}=0$ by uniqueness.

According to the definitions of $H_{jj'}^{i,*,k}, \ H_{jj'}^{i}$, we have
\begin{align*}
H_{11}^{i,*,k}(\overline P_1,\overline \Lambda_1)&\leq H_{11}^i(0,\overline P_1,\overline \Lambda_1)=0,\\
H_{12}^{i,*,k}(\overline P_1,\overline P_2,\overline \Gamma_1, \overline \Gamma_2)&\leq H_{12}^i(0,\overline P_1,\overline P_2,\overline \Gamma_1, \overline \Gamma_2),\\
H_{21}^{i,*,k}(\overline P_2,\overline \Lambda_2)&\leq H_{21}^i(0,\overline P_2,\overline \Lambda_2)=0,\\
H_{22}^{i,*,k}(\overline P_1,\overline P_2,\overline \Gamma_1, \overline \Gamma_2)&\leq H_{22}^i(0,\overline P_1,\overline P_2,\overline \Gamma_1, \overline \Gamma_2).
\end{align*}
Also, thanks to Assumption \ref{assu2},
\begin{align*}
Q^i+H_{11}^{i,*,k}(\underline P_1,\underline \Lambda_1)+\int_{\cZ}H_{12}^{i,*,k}(\underline P_1,\underline P_2,\underline \Lambda_1,\underline \Lambda_2)\nu(\dz)=Q^i\geq 0,\\
Q^i+H_{21}^{i,*,k}(\underline P_2,\underline \Lambda_2)+\int_{\cZ}H_{22}^{i,*,k}(\underline P_1,\underline P_2,\underline \Lambda_1,\underline \Lambda_2)\nu(\dz)=Q^i\geq 0.
\end{align*}
We can apply the Lemma \ref{comparison} \footnote{Conditions (1)-(5) can be obtained in a similar way as \cite[Theorem 3.1]{HSX4}} to \eqref{Ptrun} and \eqref{overlineP}, and to \eqref{Ptrun} and \eqref{underlineP}, respectively, to get
\begin{align*}
0\leq P_1^{i,k}\leq \overline P_1^{i}, \ 0\leq P_2^{i,k}\leq \overline P_2^{i}.
\end{align*}
Applying the same comparison theorem to different $k$s in \eqref{Ptrun}, we get $P_j^{i,k}$ is non-increasing in $k$, for any $i\in\cM, \ j=1,2$.

A nonnegative solution to \eqref{P} can be constructed in much the same way as \cite[Theorem 3.1]{HSX4}
by proving the strong convergence of $(P_j^{i,k},\Lambda_j^{i,k},\Gamma_j^{i,k})_{i\in\cM,\;j=1,2}$ as $k\to\infty$.
Details are left to the interested readers.

\bigskip
We now turn to the proof of uniqueness. Suppose $(P^i_j,\Lambda^i_j,\Gamma^i_j)_{i\in\cM,\;j=1,2}$ and $(\tilde P^i_j,\tilde\Lambda^i_j,\tilde\Gamma^i_j)_{i\in\cM,\;j=1,2}$ are two nonnegative solutions of \eqref{P}. Then there exists a constant $M>0$ such that
$$0\leq P^i_j,~\tilde P^i_j\leq M.$$
Estimates similar as in \cite[Theorem 3.1]{HSX4} yields also that
$$0\leq P^i_{j,t-}+\Gamma^i_{j,t},~ \tilde P^i_{j,t-}+\tilde\Gamma^j_{j,t}\leq M.$$

Firstly, we show $\int_0^{\cdot}\Lambda^i_1\dw$ is a BMO martingale.
Applying It\^o's formula to $|P_{1,t}^i|^2$, we get, for any $\mathbb{F}^{W,N}$ stopping time $\tau\leq T$,
\begin{align*}
&\quad\;\E\Big[\int_{\tau}^T|\Lambda^i_1|^2\ds\;\Big|\;\mathcal{F}^{W,N}_{\tau}\Big]\\
&\leq |G^i|^2
+\E\Big\{\int_{\tau}^T2P_{1}^i\Big[(2A^i+|C^i|^2)P^i_{1,t-}+2(C^i)^{\top}\Lambda^i_1+Q^i\\
&\qquad\qquad\qquad+H_{11}^{i,*}(P^i_1,\Lambda^i_1)+\int_{\cZ}H_{12}^{i,*}(P^i_1,P^i_2,\Gamma^i_1, \Gamma^i_2)
\nu(\dz)+\sum_{j=1}^{\ell}q^{ij}P_{1}^j\Big]\ds\;\Big|\;\mathcal{F}^{W,N}_{\tau}\Big\}\\
&\leq c+\frac{1}{2}\E\Big[\int_{\tau}^T|\Lambda^i_1|^2\ds\;\Big|\;\mathcal{F}^{W,N}_{\tau}\Big],
\end{align*}
where we used Assumption \ref{assu1}, $H_{11}^{i,*}\leq 0$, $H_{12}^{i,*}(P^i_1,P^i_2,\Gamma^i_1, \Gamma^i_2)\leq H_{12}^{i}(0,P^i_1,P^i_2,\Gamma^i_1, \Gamma^i_2)$, and the solution $(P_j^{i},\Lambda_j^{i},\Gamma_j^{i})_{i\in\cM,\;j=1,2}$ is uniformly bounded.
Note both sides in the above estimate are finite since
$\Lambda^i_1\in L^{2}_{\mathbb{F}^{W,N}}(0,T;\R^{n_1})$.
After rearrangement, we conclude that $\int_0^{\cdot}\Lambda^i_1\dw$ is a BMO martingale. Likewise, $\int_0^{\cdot}\Lambda^i_2\dw$ is also a BMO martingale.
This completes the proof of existence.

Let $a>0$ be a sufficiently small constant such that $R^i_1-a(D^i)^{\top}D^i>0$. Write $\varrho=\frac{a}{a+M}$, then $0<\varrho< 1$.
Let
\begin{align*}
&(U^i_j,V^i_j,\Phi^i_{jk})=\Big(\ln (P^i_j+a),\frac{\Lambda^i_j}{P^i_j+a},\ln\Big(1+\frac{\Gamma^i_{jk}}{P^i_{j,t-}+a}\Big)\Big),\\
&(\tilde U^i_j,\tilde V^i_j,\tilde \Phi^i_{jk})=\Big(\ln (\tilde P^i_j+a),\frac{\tilde \Lambda^i_j}{\tilde P^i_j+a},\ln\Big(1+\frac{\tilde \Gamma^i_{jk}}{\tilde P^i_{j,t-}+a}\Big)\Big)
\end{align*}
for all $i$, $j$, $k$.
Then we have the estimates
\begin{align}\label{Philowerbound}
\varrho\leq e^{\Phi^i_{jk}}, e^{\tilde\Phi^i_{jk}}\leq \varrho^{-1}, \ \ \ \frac{e^{U^i_1+\Phi^i_{1,k}}-a}{e^{U^i_1}}=\frac{P^i_1+\Gamma^i_{1k}}{e^{U^i_1}}\geq 0, \ \ \frac{e^{U^i_2+\Phi^i_{2,k}}-a}{e^{U^i_1}}=\frac{P^i_2+\Gamma^i_{2k}}{e^{U^i_1}}\geq 0.
\end{align}
Also,
\begin{align*}
\begin{cases}
\dd U_{1}^i=-\Big[(2A^i+|C^i|^2)(1-ae^{-U^i_{1}})+2(C^i)^{\top}V^i_{1}+Q^ie^{-U^i_{1}}
+\frac{1}{2}|V^i_1|^2+\sum_{j=1}^{\ell} q^{ij}e^{U_1^{j}-U_1^{i}}\\
\qquad\qquad\quad+\tilde H^{i,*}_{11}(U^i_{1},V^i_{1})+\int_{\cZ}\tilde H^{i,*}_{12}(U^i_{1},U^i_{2},\Phi^i_{1},\Phi^i_2)\nu(\dz)+\sum_{k=1}^{n_2}\int_{\cZ}(e^{\Phi^i_{1,k}}-\Phi^i_{1,k}-1)\nu(\dz)\Big]\dt\\
\qquad\quad\;+(V_{1}^i)^{\top}\dw+\int_{\cZ}\Phi_1^i(z)^{\top}\tilde N(\dt,\dz),\bigskip\\
\dd U^i_{2}=-\Big[(2A^i+|C^i|^2)(1-ae^{-U^i_{2}})+2(C^i)^{\top}V^i_{2}+Q^ie^{-U^i_{2}}
+\frac{1}{2}|V^i_2|^2+\sum_{j=1}^{\ell} q^{ij}e^{U_2^{j}-U_2^{i}}\\
\qquad\qquad\quad+\tilde H^{i,*}_{21}(U^i_{2},V^i_{2})+\int_{\cZ}\tilde H^{i,*}_{22}(U^i_{1},U^i_{2},\Phi^i_{1},\Phi^i_2)\nu(\dz)+\int_{\cZ}\sum_{k=1}^{n_2}(e^{\Phi^i_{2,k}}-\Phi^i_{2,k}-1)\nu(\dz)\Big]\dt\\
\qquad\quad\;+(V^i_{2})^{\top}\dw+\int_{\cZ}\Phi_2^i(z)^{\top}\tilde N(\dt,\dz),\bigskip\\
U^i_{1,T}=U^i_{2,T}=\ln(G^i+a),~~i\in\cM,
\end{cases}
\end{align*}
where
\begin{align*}
\tilde H^{i}_{11}(v,U_1,V_1)&:=v^{\top}(R^i_1e^{-U_{1}}+(1-ae^{-U_1})(D^i)^{\top}D^{i})v\\
&\qquad+2((1-ae^{-U_1})(B^{i}_1+(D^{i})^{\top}C^{i})+(D^{i})^{\top}V_1)^{\top}v, \\
\tilde H^{i}_{12}(v,U_1,U_2,\Phi_1, \Phi_2)&:=v^{\top}R^{i}_2e^{-U_1}v+\sum_{k=1}^{n_2}\frac{e^{U_1+\Phi_{1,k}}-a}{e^{U_1}}\Big(((1+E^{i}_k+F^{i}_kv)^+)^2-1\Big)\\
&\qquad-2(1-ae^{-U_{1}})\sum_{k=1}^{n_2}(E^{i}_k+F^{i}_kv)+2(1-ae^{-U_{1}})(B^{i}_2)^{\top}v\\
&\qquad+\sum_{k=1}^{n_2}\frac{e^{U_2+\Phi_{2,k}}-a}{e^{U_1}} ((1+E^{i}_k+F^{i}_kv)^-)^2, \\
\tilde H^{i}_{21}(v,U_2,V_2)&:=v^{\top}(R^{i}_1e^{-U_{2}}+(1-ae^{-U_2})(D^{i})^{\top}D^{i})v\\
&\qquad-2((1-ae^{-U_2})(B^{i}_1+(D^{i})^{\top}C^{i})+(D^{i})^{\top}V_2)^{\top}v, \\
\tilde H^{i}_{22}(v,U_1,U_2,\Phi_1, \Phi_2)&:=v^{\top}R^{i}_2e^{-U_2}v+\sum_{k=1}^{n_2}\frac{e^{U_2+\Phi_{2,k}}-a}{e^{U_2}}\Big(((-1-E^{i}_k+F^{i}_kv)^-)^2-1\Big)\\
&\qquad-2(1-ae^{-U_{2}})\sum_{k=1}^{n_2}(-E^{i}_k+F^{i}_kv)-2(1-ae^{-U_{2}})(B_2^{i})^{\top}v\\
&\qquad+\sum_{k=1}^{n_2}\frac{e^{U_1+\Phi_{1,k}}-a}{e^{U_2}} ((-1-E^{i}_k+F^{i}_kv)^-)^2,
\end{align*}
and
\begin{align*}
\tilde H_{11}^{i,*}(U_1,V_1)&:=\inf_{v\in\Pi_1}
\tilde H^{i}_{11}(v,U_1,V_1),\\
\tilde H_{12}^{i,*}(z,U_1,U_2,\Phi_1, \Phi_2)&:=\inf_{v\in\Pi_2}
\tilde H^{i}_{12}(v,z,U_1,U_2,\Phi_1, \Phi_2),\\
\tilde H_{21}^{i,*}(U_2,V_2)&:=\inf_{v\in\Pi_1}
\tilde H^{i}_{21}(v,U_2,V_2),\\
\tilde H_{22}^{i,*}(z,U_1,U_2,\Phi_1, \Phi_2)&:=\inf_{v\in\Pi_2}
\tilde H^{i}_{22}(v,z,U_1,U_2,\Phi_1, \Phi_2).
\end{align*}

Set
\[
\bar U^i_j=U^i_j-\tilde U^i_j, \ \bar V^i_j=V^i_j-\tilde V^i_j, \ \bar\Phi^i_j=\Phi^i_j-\tilde\Phi^i_j, \ i\in\cM, \ j=1,2.
\]
Then applying It\^{o}'s formula to $(\bar U^i_j)^2$, we deduce that
\begin{align*}
&\quad\;(U^i_{1,t})^2+\int_t^T|\bar V^i_1|^2\ds+\int_t^T\int_{\cZ}|\bar\Phi^i_1|^2\nu(\dz)\dt\\
&=\int_t^T\Big[L^i_{11}+\int_{\cZ}L^i_{12}(z)\nu(\dz)\Big]\dt-\int_t^T2\bar U^i_1(\bar V^i_{1})^{\top}\dw\\
&\quad\;-\sum_{k=1}^{n_2}\int_t^T\int_{\cZ}(2\bar U^i_{1,k}\bar\Phi^i_{1,k}+(\bar\Phi^i_{1,k})^2)\tilde N_k(\dt,\dz),
\end{align*}
and
\begin{align*}
&\quad\;(U^i_{2,t})^2+\int_t^T|\bar V^i_2|^2ds+\int_t^T\int_{\cZ}|\bar\Phi^i_2|^2\nu(\dz)\dt\\
&=\int_t^T\Big[L^i_{21}+\int_{\cZ}L^i_{22}(z)\nu(\dz)\Big]\dt
-\int_t^T2\bar U^i_2(\bar V^i_{2})^{\top}\dw\\
&\quad\;-\sum_{k=1}^{n_2}\int_t^T\int_{\cZ}(2\bar U^i_{2,k}\bar\Phi^i_{2,k}+(\bar\Phi^i_{2,k})^2)\tilde N_k(\dt,\dz),
\end{align*}
where
\begin{align*}
L^i_{11}:&=2\bar U^i_1\Big[(Q^i-2aA^i-a|C^i|^2)(e^{-U^i_1}-e^{-\tilde U^i_1})+2(C^i)^{\top}\bar V^i_1+\frac{1}{2}(V^i_1+\tilde V^i_1)\bar V^i_1\\
&\qquad\qquad+\sum_{j=1}^{\ell}q^{ij}\Big(e^{U^i_1-U^i_1}-e^{\tilde U^j_1-\tilde U^i_1}\Big)+\tilde H_{11}^{i,*}(U^i_{1},V^i_{1},)-\tilde H_{11}^{i,*}(\tilde U^i_{1},\tilde V^i_{1})\Big],\\
L^i_{12}(z):&=2\bar U^i_1\Big[\sum_{k=1}^{n_2}[(e^{\Phi^i_{1,k}}-\Phi^i_{1,k}-1)-(e^{\tilde\Phi^i_{1,k}}-\tilde\Phi^i_{1,k}-1)]\\
&\qquad\qquad +\tilde H_{12}^{i,*}(z,U^i_{1},\Phi^i_1,U^i_{2},\Phi^i_2)-\tilde H_{12}^{i,*}(z,\tilde U^i_{1},\tilde\Phi^i_1,\tilde U^i_{2},\tilde\Phi^i_2)\Big],\\
L^i_{21}:&=2\bar U^i_2\Big[(Q^i-2aA^i-a|C^i|^2)(e^{-U^i_2}-e^{-\tilde U^i_2})+2(C^i)^{\top}\bar V^i_2+\frac{1}{2}(V^i_2+\tilde V^i_2)\bar V^i_2\\
&\qquad\qquad+\sum_{j=1}^{\ell}q^{ij}\Big(e^{U^i_2-U^i_2}-e^{\tilde U^j_2-\tilde U^i_2}\Big)+\tilde H_{21}^{i,*}(U^i_{2},V^i_{2},)-\tilde H_{21}^{i,*}(\tilde U^i_{2},\tilde V^i_{2})\Big],\\
L^i_{22}(z):&=2\bar U^i_2\Big[\sum_{k=1}^{n_2}[(e^{\Phi^i_{2,k}}-\Phi^i_{2,k}-1)-(e^{\tilde\Phi^i_{2,k}}-\tilde\Phi^i_{2,k}-1)]\\ &\qquad\qquad+\tilde H_{22}^{i,*}(z,U^i_{1},\Phi^i_1,U^i_{2},\Phi^i_2)-\tilde H_{22}^{i,*}(z,\tilde U^i_{1},\tilde\Phi^i_1,\tilde U^i_{2},\tilde\Phi^i_2)\Big].
\end{align*}

The terms $L_{11}$ and $L_{21}$ can be estimated in much the same way as \cite[Theorem 3.5]{HSX} to get\footnote{$R^i_1-a(D^i)^{\top}D^i>0$ is required in these estimates.}
\begin{align*}
&L^i_{11}\leq |\beta^i|(\bar U^i_1)^2+c|\bar U^i_1|\sum_{j=1}^{\ell}|\bar U^j_1|+c(\beta^i)^{\top}\bar U^i_1\bar V^i_1,\\
&L^i_{21}\leq |\beta^i|(\bar U^i_2)^2+c|\bar U^i_2|\sum_{j=1}^{\ell}|\bar U^j_2|+c(\beta^i)^{\top}\bar U^i_2\bar V^i_2,
\end{align*}
where $\beta^i$ is some $\mathbb{F}^{W,N}$-predictable process satisfying $|\beta^i|\leq c(1+|V^i_1|+|\tilde V^i_1|+|V^i_2|+|\tilde V^i_2|)$ so that $\int_0^{\cdot}(\beta^i)^{\top}\dw$ is a BMO martingale.

On the other hand, from Assumptions \ref{assu1}, \ref{assu2} and \eqref{Philowerbound}, there are positive constants $c_1,c_2,c_3$ such that
\begin{align*}
&\qquad\tilde H^i_{12}(v,z,U^i_1,U^i_2,\Phi^i_1, \Phi^i_2)-\tilde H^i_{12}(0,z,U^i_1,U^i_2,\Phi^i_1, \Phi^i_2)\\
&\geq \frac{\delta}{M+a}|v|^2-\sum_{k=1}^{n_2}\frac{e^{U^i_1+\Phi^i_{1,k}}-a}{e^{U_1}} -2(1-ae^{-U^i_{1}})\sum_{k=1}^{n_2}(E^i_k+F^i_kv)+2(1-ae^{-U^i_{1}})(B^i_2)^{\top}v\\
&\qquad-\Big[\sum_{k=1}^{n_2}\frac{e^{U^i_1+\Phi^i_{1,k}}-a}{e^{U_1}}\Big(((1+E^i_k)^+)^2-1\Big)\\
&\qquad\qquad-2(1-ae^{-U^i_{1}})\sum_{k=1}^{n_2}E^i_k+\sum_{k=1}^{n_2}\frac{e^{U^i_2+\Phi^i_{2,k}}-a}{e^{U^i_1}} ((1+E^i_k)^-)^2\Big]\\
&\geq \frac{\delta}{M+a}|v|^2-c_2|v|-c_3> 0,
\end{align*}
if $|v|> c$ with $c>0$ being sufficiently large. Hence,
\begin{align}\label{boundeddomain}
\tilde H_{12}^{i,*}(z,U^i_1,U^i_2,\Phi^i_1, \Phi^i_2)&:=\inf_{v\in\Pi_2, |v|\leq c}
\tilde H^i_{12}(v,z,U^i_1,U^i_2,\Phi^i_1, \Phi^i_2).
\end{align}
Furthermore, noting $U^i_1,\tilde U^i_1, U^i_2,\tilde U^i_2, \Phi^i_1, \tilde\Phi^i_1, \Phi^i_2, \tilde\Phi^i_2$ are bounded, we have
\begin{align*}
L^i_{12}(z) \leq c|\bar U^i_1|(|\bar U^i_1|+|\bar\Phi^i_1(z)|+|\bar U^i_2|+|\bar\Phi^i_2(z)|).
\end{align*}
Similar arguments applying to $\tilde H^{i,*}_{22}(v,z,U^i_1,U^i_2,\Phi^i_1, \Phi^i_2)$ yield that
\begin{align*}
L^i_{22}(z) \leq c|\bar U^i_2|(|\bar U^i_1|+|\bar\Phi^i_1(z)|+|\bar U^i_2|+|\bar\Phi^i_2(z)|).
\end{align*}

For each $i\in\cM$, introduce the processes
\begin{align*}
J^i_t=\exp\Big(\int_0^t|\beta^i_s|\ds\Big), \ N^i_t=\exp\Big(\int_0^t(\beta^i_s)^{\top}\dw_s-\frac{1}{2}\int_0^t|\beta^i_s|^2\ds\Big).
\end{align*}
It\^{o}'s formula gives
\begin{align*}
&\qquad J^i_tN^i_t|\bar U^i_{1,t}|^2+\E_t\int_t^TJ^i_sN^i_s|\bar V^i_1|^2\ds+\E_t\int_t^T\int_{\cZ}J^i_sN^i_s|\bar\Phi^i_1|^2\nu(\dz)\ds\\
&\leq\E_t\int_t^TJ^i_sN^i_s\Big[c|\bar U^i_1||\bar U^i_2|+c|\bar U^i_1|\sum_{j=1}^{\ell}|\bar U^j_1|+c|U^i_1|\int_{\cZ}(\bar\Phi^i_1(z)+\bar\Phi^i_2(z))\nu(\dz)\Big]\ds\\
&\leq c\E_t\int_t^{T}J^i_sN^i_s(|\bar U^i_1|^2+|\bar U^i_2|^2+\sum_{j=1}^{\ell}|\bar U_1^j|^2)\ds+\frac{1}{4}\E_t\int_t^T\int_{\cZ}J^i_sN^i_s(|\bar\Phi^i_1(z)|^2+|\bar\Phi^i_2(z)|^2)\nu(\dz)\ds.
\end{align*}
Note that $N^i_t$ is a uniformly integrable martingale, thus
\begin{align*}
\widetilde W^i_t:=W_t-\int_0^t(\beta^i_s)^{\top}\dw_s,
\end{align*}
is a Brownian motion under the probability $\tilde{\mathbb{P}}^i$ defined by
\[
\frac{d\tilde{\mathbb{P}}^i}{d\mathbb{P}}\Big|_{\mathcal{F}^{W,N}_T}=N^i_T.
\]
We denote by $\widetilde \E^i_t$ the conditional expectation with respect to the probability $\tilde{\mathbb{P}}^i$, then
\begin{align}\label{e1}
&\qquad J^i_t|\bar U^i_{1,t}|^2+\widetilde\E^i_t\int_t^TJ^i_s|\bar V^i_1|^2\ds+\widetilde\E^i_t\int_t^T\int_{\cZ}J^i_s|\bar\Phi^i_1|^2\nu(\dz)\ds\nn\\
&\leq c\widetilde\E^i_t\int_t^{T}J^i_s(|\bar U^i_1|^2+|\bar U^i_2|^2+\sum_{j=1}^{\ell}|\bar U_1^j|^2)\ds+\frac{1}{4}\widetilde\E^i_t\int_t^T\int_{\cZ}J^i_s(|\bar\Phi^i_1|^2+|\bar\Phi^i_2|^2)\nu(\dz)\ds.
\end{align}
Similarly, we have
\begin{align}\label{e2}
&\qquad J^i_t|\bar U^i_{2,t}|^2+\widetilde\E^i_t\int_t^TJ^i_s|\bar V^i_2|^2\ds+\widetilde\E^i_t\int_t^T\int_{\cZ}J^i_s|\bar\Phi^i_2|^2\nu(\dz)\ds\nn\\
&\leq c\widetilde\E^i_t\int_t^{T}J^i_s(|\bar U^i_1|^2+|\bar U^i_2|^2+\sum_{j=1}^{\ell}|\bar U_2^j|^2)\ds+\frac{1}{4}\widetilde\E^i_t\int_t^T\int_{\cZ}J^i_s(|\bar\Phi^i_1|^2+|\bar\Phi^i_2|^2)\nu(\dz)\ds.
\end{align}
Combining the above two inequalities yields
\begin{align*}
|\bar U_{1,t}^i|^2+|\bar U^i_{2,t}|^2&\leq c\widetilde\E^i_t\int_t^{T}\exp(\int_t^s|\beta^i_r|\dr)(\sum_{j=1}^{\ell}\int_t^{T}(|\bar U^j_{1,s}|^2+|\bar U^j_{2,s}|^2)\ds\\
&\leq c\widetilde\E^i_t \Big[\exp(\int_t^T|\beta^i_r|\dr) \sum_{j=1}^{\ell}\int_t^{T}(|\bar U^j_{1,s}|^2+|\bar U^j_{2,s}|^2)\ds\Big]\\
&\leq c\widetilde\E^i_t \Big[\exp(\int_t^T|\beta^i_r|\dr)\Big]\sum_{j=1}^{\ell}\int_t^{T}\Xi^j_{s}\ds,
\end{align*}
where
\[
\Xi^j_{s}:=\esssup_{\omega\in\Omega} \Big(|\bar U^j_{1,s}|^{2}+|\bar U^j_{2,s}|^{2}\Big).
\]
According to \cite[Lemma 3.4]{HSX3}, $\widetilde\E^i_t \Big[\exp(\int_t^T|\beta^i_r|\dr)\Big]\leq c$, then taking essential supreme on both sides, we deduce that
\begin{align*}
0\leq \sum_{i=1}^{\ell}\Xi^i_{t}\leq c\int_t^{T}\sum_{i=1}^{\ell}\Xi^i_{s} \ds.
\end{align*}
We infer from Gronwall's inequality that $\Xi^i=0$, so $\bar U^i_{1}=\bar U^i_{2}=0$, for all $i\in\cM$. Consequently, it follows from \eqref{e1} and \eqref{e2} that $\bar V^i_{1}=\bar V^i_{2}=0$ and $\bar\Phi^i_1=\bar\Phi^i_2=0$ for all $i\in\cM$. This completes the proof.
\eof

\begin{theorem}\label{existencesing}
Under Assumptions \ref{assu1} and \ref{assu3}, the BSDEJ \eqref{P} admits a unique uniformly positive solution $(P^i_j,\Lambda^i_j,\Gamma^i_j)_{i\in\cM,\;j=1,2}$.
\end{theorem}
\pf
The proof of the existence is similar to the above Theorem \ref{existence} and will only be indicated briefly why the solution to \eqref{P} is uniformly positive.

\underline{When both \ref{R12} and \ref{R21} hold.}
In this case, there exists constant $c_2>0$, such that
\[
2A^i+|C^{i}|^2-\delta^{-1}|B^i_1+(D^i)^{\top}C^i|^2\geq -c_2, \ -\delta^{-1}|B_2^i|^2\geq -c_2
\]
where $\delta$ is the constant in Assumption \ref{assu3}.
Notice that $(\underline P_{t},\underline\Lambda_{t},\underline\Gamma_{t})=(\frac{1}{(\delta^{-1}+1)e^{c_2(T-t)}-1},0,0)$ solves the following BSDEJ
\begin{align}\label{Plowersingu1}
\begin{cases}
\dd\underline P=-(-c_2\underline P-c_2\underline P^2)\dt+\underline\Lambda^{\top}\dw+\int_{\cE}\underline\Gamma(e) \widetilde N(\dt,\dz),\\
\underline P_{T}=\delta.
\end{cases}
\end{align}
And we have the following inequalities
\begin{align*}
H_{11}^{i,*,k}(\underline P,\underline \Lambda)&\geq \inf_{v\in\R^{m_1}}
H_{11}^{i}(v,\underline P,\underline \Lambda)\geq-\delta^{-1} |B^i_1+(D^i)^{\top}C^i|^2\underline P,\\
H_{12}^{i,*,k}(\underline P,\underline P,\underline \Gamma,\underline \Gamma)&\geq \inf_{v\in\R^{m_1}}H_{12}^{i}(v,\underline P,\underline P,\underline \Gamma,\underline \Gamma)\geq-\delta^{-1}|B_2^i|^2\underline P^2\geq -c_2\underline P^2.
\end{align*}
Similar results also holds for $H_{21}^{i,*,k}(\underline P,\underline \Lambda)$ and $H_{22}^{i,*,k}(\underline P,\underline P,\underline \Gamma,\underline \Gamma)$.
Applying Lemma \ref{comparison} to \eqref{Ptrun} and \eqref{Plowersingu1}, we get
\begin{align*}
P^{i,k}_{j,t}\geq \underline P_{t}\geq \frac{1}{(\delta^{-1}+1)e^{c_2T}-1}, \ t\in[0,T], \ i\in\cM, \ j=1,2.
\end{align*}
Sending $k\rightarrow\infty$ leads to the desired uniformly positive lower bound.

\underline{When both \ref{R12} and \ref{R22} hold.}
In this case, there exists constant $c_3>0$, such that
\[
2A^i+|C^{i}|^2-\delta^{-1}|B^i_1+(D^i)^{\top}C^i|^2-\delta^{-1}\int_{\cZ}|(F^i)^{\top}E^i+B_2^i|^2\nu(\dz)\geq -c_3,
\]
where $\delta$ is the constant in Assumption \ref{assu3}.
Notice $(\underline P_{t},\underline\Lambda_{t},\underline\Gamma_{t})=(\delta e^{-c_3(T-t)},0,0)$ solves the following BSDEJ
\begin{align}\label{Plowersingu2}
\begin{cases}
\dd\underline P=-(-c_3\underline P)\dt+\underline\Lambda^{\top}\dw+\int_{\cE}\underline\Gamma(e) \widetilde N(\dt,\dz),\\
\underline P_{T}=\delta.
\end{cases}
\end{align}
And we have the following inequalities
\begin{align*}
H_{11}^{i,*,k}(\underline P,\underline \Lambda)&\geq \inf_{v\in\R^{m_1}}
H_{11}^{i}(v,\underline P,\underline \Lambda)\geq-\delta^{-1} |B^i_1+(D^i)^{\top}C^i|^2\underline P,\\
H_{12}^{i,*,k}(z,\underline P,\underline P,\underline \Gamma,\underline \Gamma)&\geq \inf_{v\in\R^{m_1}}H_{12}^{i}(z, v,\underline P,\underline P,\underline \Gamma,\underline \Gamma)\geq-\delta^{-1}|(F^i)^{\top}E^i+B_2^i|^2\underline P.
\end{align*}
Therefore,
\begin{align*}
&(2A^i+|C^i|^2)\underline P+2(C^i)^{\top}\underline\Lambda+Q^i+\sum_{j=1}^{\ell}q^{ij}\underline P+H_{11}^{i,*,k}(\underline P,\underline \Lambda)+\int_{\cZ}H_{12}^{i,*,k}(\underline P,\underline P,\underline \Gamma,\underline \Gamma)\nu(\dz)\\
&\geq (2A^i+|C^i|^2)\underline P-\delta^{-1} |B^i_1+(D^i)^{\top}C^i|^2\underline P-\int_{\cZ}\delta^{-1}|(F^i)^{\top}E^i+B_2^i|^2\underline P\nu(\dz)\geq -c_3\underline P.
\end{align*}

Similar results also holds for $H_{21}^{i,*,k}(\underline P,\underline \Lambda)$ and $H_{22}^{i,*,k}(z,\underline P,\underline P,\underline \Gamma,\underline \Gamma)$.
Applying Lemma \ref{comparison} to \eqref{Ptrun} and \eqref{Plowersingu2}, we get
\begin{align*}
P^{i,k}_{j,t}\geq \underline P_{t}=\delta e^{-c_3(T-t)}\geq \delta e^{-c_3T}, \ t\in[0,T], \ i\in\cM, \ j=1,2.
\end{align*}
Sending $k\rightarrow\infty$ leads to the desired uniformly positive lower bound.

\underline{When both \ref{R11} and \ref{R22} hold.}
This case can be handled in a similar manner. The details are left to the readers.

As for the uniqueness, just setting $a=0$ in the proof of the uniqueness part of Theorem \ref{existence}.
\eof

\section{Solution to the LQ problem \eqref{LQ}}\label{section:Veri}
In this subsection we will present an explicit solution to the LQ problem \eqref{LQ} in terms of solutions to the BSDEJ \eqref{P}.

When $R^i_1+P(D^i)^{\top}D^i>0, \ i\in\cM$, we define
\begin{align}\label{hatv1}
\hat v_{11}^{i}(t,P,\Lambda)&:=\argmin_{v\in\Pi_1}
H_{11}^i(t,v,P,\Lambda),\nn\\
\hat v_{21}^{i}(t,P,\Lambda)&:=\argmin_{v\in\Pi_1}
H_{21}^i(t,v,P,\Lambda).
\end{align}
When $R^i_2>0$ or $(P_j+\Gamma_{jk})(F^i)^{\top}F^i>0$, $i\in\cM, \ j=1,2, \ k=1,2,...,n_2$, we
define
\begin{align}\label{hatv2}
\hat v_{12}^{i}(t,z,P_1,P_2,\Gamma_1, \Gamma_2)&:=\argmin_{v\in\Pi_2}
H_{12}^i(t,v,z,P_1,P_2,\Gamma_1, \Gamma_2),\nn\\
\hat v_{22}^{i}(t,z,P_1,P_2,\Gamma_1, \Gamma_2)&:=\argmin_{v\in\Pi_2}
H_{22}^i(t,v,z,P_1,P_2,\Gamma_1, \Gamma_2).
\end{align}


\begin{theorem}\label{Th:verif}
Let Assumptions \ref{assu1}, and \ref{assu2} (resp. \ref{assu3}) hold.
Let $(P^i_j,\Lambda^i_j,\Gamma^i_j)\in S^{\infty}_{\mathbb{F}^{W,N}}(0,T;\R)\times L^{2}_{\mathbb{F}^{W,N}}(0,T;\R^{n_1}) \times L^{\infty,\nu}_{\mathcal{P}^{W,N}}(0, T;\R^{n_2}), \ i\in\cM, \ j=1,2$, be the nonnegative (resp. uniformly positive) solution to the BSDEJ \eqref{P}. Then the state feedback control $u^*=(u^*_1,u^*_2)$ given by
\begin{align} \label{ustar}
\begin{cases}
u^*_1(t,X,\alpha)=\hat v_{11}^{\alpha_{t-}}(t,P^{\alpha_{t-}}_{1,t-},\Lambda^{\alpha_{t-}}_{1,t})X_{t-}^++\hat v_{21}^{\alpha_{t-}}(t,P^{\alpha_{t-}}_{2,t-},\Lambda^{\alpha_{t-}}_{2,t})X_{t-}^-,\\
u^*_2(t,X,\alpha)=\hat v_{12}^{\alpha_{t-}}(t,z,P^{\alpha_{t-}}_{1,t-},P^{\alpha_{t-}}_{2,t-},\Gamma^{\alpha_{t-}}_{1,t}, \Gamma^{\alpha_{t-}}_{2,t})X_{t-}^+ +\hat v_{22}^{\alpha_{t-}}(t,z,P^{\alpha_{t-}}_{1,t-},P^{\alpha_{t-}}_{2,t-},\Gamma^{\alpha_{t-}}_{1,t}, \Gamma^{\alpha_{t-}}_{2,t})X_{t-}^-,
\end{cases}
\end{align}
is optimal for the LQ problem \eqref{LQ}. Moreover, the optimal value is
\begin{align*}
V(x,i_0)=P^{i_0}_{1,0}(x^+)^2+P^{i_0}_{2,0}(x^-)^2.
\end{align*}
\end{theorem}

A proof of this theorem is contained in the following two Lemmas.
In order to avoid as far as possible unwieldy formulas, we agree to suppress the the superscripts and subscripts of $A,B,C,D,E,F,R,Q,G$. And we will write $v_{ij}^{\alpha_{t-}}(t,P^{\alpha_{t-}}_1,\Lambda^{\alpha_{t-}}_1)$ simply $\hat v_{ij}$, $i,j=1,2$ when no confusion can arise.

\begin{lemma}\label{lemma1}
Under the condition of Theorem \ref{Th:verif}, we have
\begin{align}\label{lower1}
J(u;x,i_0)\geq P^{i_0}_{1,0}(x^+)^2+P^{i_0}_{2,0}(x^-)^2,
\end{align}
for any $u\in\mathcal{U}$, and
\begin{align}
J(u^*;x,i_0)= P^{i_0}_{1,0}(x^+)^2+P^{i_0}_{2,0}(x^-)^2.
\end{align}
\end{lemma}
\pf
By the Meyer-It\^o formula \cite[Theorem 70]{Protter}, we deduce that, for any $u=(u_1,u_2)\in\mathcal{U}$,
\begin{align*}
\dd X_t^+&=\mathbf{1}_{\{X_{t-}>0\}}\Big[\Big(AX+B_1^{\top}u_1
+\int_{\cZ}B_2(z)^{\top}u_2(z)\nu(\dz)-\sum_{k=1}^{n_2}\int_{\cZ}(E_k(z)X+F_k(z)u_2(z))\nu(\dz)\Big)\dt\\
&\quad+(CX+Du_1)^{\top}\dw_t\Big]+\sum_{k=1}^{n_2}\int_{\cZ}\Big[((X+E_k(z)X+F_k(e)u_2(z))^+-X^+\Big]N_k(\dt,\dz)+\frac{1}{2}\dd L_t,
\end{align*}
and
\begin{align*}
\dd X_t^-&=-\mathbf{1}_{\{X_{t-}\leq0\}}\Big[\Big(AX+B_1^{\top}u_1
+\int_{\cZ}B_2(z)^{\top}u_2(z)\nu(\dz)-\sum_{k=1}^{n_2}\int_{\cZ}(E_k(z)X+F_k(z)u_2(z))\nu(\dz)\Big)\dt\\
&\quad+(CX+Du_1)^{\top}\dw_t\Big]+\sum_{k=1}^{n_2}\int_{\cZ}\Big[((X+E_k(z)X+F_k(e)u_2(z))^--X^-\Big]N_k(\dt,\dz)+\frac{1}{2}\dd L_t,
\end{align*}
where $L$ is the local time of $X$ at $0$. Since $X_{t-}^{\pm}\dd L_t=0$, it follows from the It\^{o} formula that
\begin{align*}
\dd\; (X_t^+)^2&=\Big[2X_{t-}^+\Big(AX+B_1^{\top}u_1
+\int_{\cZ}B_2(z)^{\top}u_2(z)\nu(\dz)-\sum_{k=1}^{n_2}\int_{\cZ}(E_k(z)X+F_k(z)u_2(z))^{\top}\nu(\dz)\Big)\\
&\qquad+\mathbf{1}_{\{X_{t-}>0\}}|CX+Du_1|^2\Big]\dt+2X^+(CX+Du_1)^{\top}\dw_t\\
&\qquad+\sum_{k=1}^{n_2}\int_{\cZ}\Big[((X+E_k(z)X+F_k(e)u_2(z))^+)^2-(X^+)^2\Big]N_k(\dt,\dz),
\end{align*}
and
\begin{align*}
\dd\;(X_t^-)^2&=\Big[-2X_{t-}^-\Big(AX+B_1^{\top}u_1
+\int_{\cZ}B_2(z)^{\top}u_2(z)\nu(\dz)-\sum_{k=1}^{n_2}\int_{\cZ}(E_k(z)X+F_k(z)u_2(z))^{\top}\nu(\dz)\Big)\\
&\qquad+\mathbf{1}_{\{X_{t-}\leq0\}}|CX+Du_1|^2\Big]\dt-2X^-(CX+Du_1)^{\top}\dw_t\\
&\qquad+\sum_{k=1}^{n_2}\int_{\cZ}\Big[((X+E_k(z)X+F_k(e)u_2(z))^-)^2-(X^-)^2\Big]N_k(\dt,\dz).
\end{align*}
Then applying the It\^{o} formula to $P^{\alpha_t}_{1,t}(X_t^+)^2$, yields
\begin{align*}
\dd P^{\alpha_t}_{1,t}(X_t^+)^2&=P^{\alpha_{t-}}_{1,t-}\Big[2X_{t-}^+\Big(AX+B_1^{\top}u_1+\int_{\cZ}B_2(z)^{\top}u_2\nu(\dz)\\
&\qquad-\sum_{k=1}^{n_2}\int_{\cZ}(E_k(z)X+F_k(z)u_2(z))\nu(\dz)\Big)
+\mathbf{1}_{\{X_{t-}>0\}}|CX+Du_1|^2\Big]\dt\\
&\qquad+2X^+(CX+Du_1)^{\top}\Lambda^{\alpha_{t-}}_1\dt\\
&\qquad-(X^+)^2\Big[(2A+|C|^2)P^{\alpha_{t-}}_{1,t-}+2C^{\top}\Lambda^{\alpha_{t-}}_1+Q
+H_{11}^{\alpha_{t-},*}(P^{\alpha_{t-}}_1,\Lambda^{\alpha_{t-}}_1)\\
&\qquad\qquad+\int_{\cZ}H_{12}^{\alpha_{t-},*}(P^{\alpha_{t-}}_1,P^{\alpha_{t-}}_2,\Gamma^{\alpha_{t-}}_1, \Gamma^{\alpha_{t-}}_2)\nu(\dz)\Big]\dt\\
&\qquad+\sum_{k=1}^{n_2}\int_{\cZ}(P^{\alpha_{t-}}_{1,t-}+\Gamma^{\alpha_{t-}}_{1k,t})\Big[\big((X+E_k(z)X+F_k(e)u_2(z))^+\big)^2-(X^+)^2\Big]\nu(\dz)\dt\\
&\qquad +\Big[2X^+(CX+Du_1)+(X^+)^2\Lambda^{\alpha_{t-}}_1\Big]^{\top}\dw\\
&\qquad+\sum_{k=1}^{n_2}\int_{\cZ}(P^{\alpha_{t-}}_{1,t-}+\Gamma^{\alpha_{t-}}_{1k,t})\Big[\big((X+E_k(z)X+F_k(e)u_2(z))^+\big)^2-(X^+)^2\Big]\tilde N_k(\dt,\de)\\
&\qquad + (X^+)^2\int_{\cZ}\Gamma^{\alpha_{t-}}_1(z)^{\top}\tilde N(\dt,\dz)+(X^+)^2 \sum_{j,j'\in\cM}(P_1^{j}-P_1^{j'})\mathbf{1}_{\{\alpha_t-=j'\}}\dd \tilde N^{j'j},
\end{align*}
where $\{N^{j'j}\}_{j,j'\in\cM}$ are independent Poisson processes each with intensity $q^{j'j}$, and $\tilde N^{j'j}_t=N^{j'j}_t-q^{j'j}t, \ t\geq 0$ are the corresponding compensated Poisson martingale.
Likewise,
\begin{align*}
\dd P^{\alpha_t}_{2,t}(X_t^-)^2&=P^{\alpha_{t-}}_{2,t-}\Big[-2X_{t-}^+\Big(AX+B_1^{\top}u_1+\int_{\cZ}B_2(z)^{\top}u_2\nu(\dz)\\
&\qquad-\sum_{k=1}^{n_2}\int_{\cZ}(E_k(z)X+F_k(z)u_2(z))\nu(\dz)\Big)
+\mathbf{1}_{\{X_{t-}\leq0\}}|CX+Du_1|^2\Big]\dt\\
&\qquad-2X^-(CX+Du_1)^{\top}\Lambda^{\alpha_{t-}}_1\dt\\
&\qquad-(X^-)^2\Big[(2A+|C|^2)P^{\alpha_{t-}}_{2,t-}+2C^{\top}\Lambda^{\alpha_{t-}}_2+Q
+H_{21}^{\alpha_{t-},*}(P^{\alpha_{t-}}_2,\Lambda^{\alpha_{t-}}_2)\\
&\qquad\qquad+\int_{\cZ}H_{22}^{\alpha_{t-},*}(P^{\alpha_{t-}}_1,P^{\alpha_{t-}}_2,\Gamma^{\alpha_{t-}}_1, \Gamma^{\alpha_{t-}}_2)\nu(\dz)\Big]\dt\\
&\qquad+\sum_{k=1}^{n_2}\int_{\cZ}(P^{\alpha_{t-}}_{2,t-}+\Gamma^{\alpha_{t-}}_{2k,t})\Big[\big((X+E_k(z)X+F_k(e)u_2(z))^-\big)^2-(X^-)^2\Big]\nu(\dz)\dt\\
&\qquad +\Big[-2X^-(CX+Du_1)+(X^-)^2\Lambda^{\alpha_{t-}}_1\Big]^{\top}\dw\\
&\qquad+\sum_{k=1}^{n_2}\int_{\cZ}(P^{\alpha_{t-}}_{2,t-}+\Gamma^{\alpha_{t-}}_{2k,t})\Big[\big((X+E_k(z)X+F_k(e)u_2(z))^-\big)^2-(X^-)^2\Big]\tilde N_k(\dt,\dz)\\
&\qquad + (X^-)^2\int_{\cZ}\Gamma^{\alpha_{t-}}_2(z)^{\top}\tilde N(\dt,\dz)+(X^-)^2 \sum_{j,j'\in\cM}(P_2^{j}-P_2^{j'})\mathbf{1}_{\{\alpha_t-=j'\}}\dd \tilde N^{j'j}.
\end{align*}
We define, for $n\geq 1$, the following stopping time $\tau_n$:
\begin{align*}
\tau_n:=\inf\{t\geq 0 : |X_{t}|\geq n\}\wedge T,
\end{align*}
with the convention that $\inf\emptyset=\infty$. Obviously, $\tau_n\uparrow T$ a.s. along $n\uparrow\infty$.

Summing the two equations above, takeing integration from $0$ to $\tau_n$, and then taking expectation, we deduce
\begin{align}\label{lower2}
&\E\Big[ P^{\alpha_{\tau_n}}_{1,{\tau_n}}(X_{\tau_n}^+)^2+P^{\alpha_{\tau_n}}_{2,{\tau_n}}(X_{\tau_n}^-)^2\Big]\nn\\
&\qquad\qquad+\E\int_0^{\tau_n}\Big[u_{1}^{\top}R_{1}u_{1}
+QX^2+\int_{\cZ}(u_{2}(z)^{\top}R_{2}(z)u_{2}(z))\nu(\dz)\Big]\dt\nn\\
&=P^{i_0}_{1,0}(x^+)^2+P^{i_0}_{2,0}(x^-)^2
+\E\int_0^{\tau_n}\Big[u_1^{\top}(R_1+\mathbf{1}_{\{X>0\}}P_1D^{\top}D+\mathbf{1}_{\{X\leq0\}}P_2D^{\top}D)u_1\nn\\
&\qquad+2u_1^{\top}(P_1B_1+D^{\top}(P_1C+\Lambda_1))X^+
-2u_1^{\top}(P_2B_1+D^{\top}(P_2C+\Lambda_2))X^-\nn\\
&\qquad -H_{11}^{\alpha_{t-},*}(P_1,\Lambda_1)(X^+)^2-H_{21}^{\alpha_{t-},*}(P_2,\Lambda_2)(X^-)^2\Big]\dt\nn\\
&\qquad + \E\int_0^{\tau_n}\int_{\cZ}\Big[u_2^{\top}R_2u_2 +2P_1X^+\Big(B_2(z)^{\top}u_2(z)-\sum_{k=1}^{n_2}(E_k(z)X+F_k(z)u_2(z))\Big)\nn\\ &\qquad\qquad-2P_2X^-\Big(B_2(z)^{\top}u_2(z)-\sum_{k=1}^{n_2}(E_k(z)X+F_k(z)u_2(z))\Big)\nn\\
&\qquad\qquad+\sum_{k=1}^{n_2}(P_{1}+\Gamma_{1k})\Big(\big((X+E_k(z)X+F_k(e)u_2(z))^+\big)^2-(X^+)^2\Big)\nn\\
&\qquad\qquad+\sum_{k=1}^{n_2}(P_{2}+\Gamma_{2k})\Big(\big((X+E_k(z)X+F_k(e)u_2(z))^-\big)^2-(X^-)^2\Big)\nn\\
&\qquad\qquad-H_{12}^{\alpha_{t-},*}(P_1,P_2,\Gamma_1, \Gamma_2)(X^+)^2
-H_{22}^{\alpha_{t-},*}(P_1,P_2,\Gamma_1, \Gamma_2)(X^-)^2\Big]\nu(\dz)\dt.
\end{align}
We will denote by $\phi(X,u)$ the right-hand side (RHS) of the above equation and show $\phi(X,u)\geq 0, ~\dptv$, for any $u\in\mathcal{U}$.

Indeed, let us define
\begin{align*}
v_t=(v_{1,t},v_{2,t}(z))=
\begin{cases}
\Big(\frac{u_{1,t}}{|X_{t-}|},\frac{u_{2,t}(z)}{|X_{t-}|}\Big), & \ \mbox {if} \ |X_{t-}|>0;\\
(0,0), & \ \mbox {if} \ |X_{t-}|=0.
\end{cases}
\end{align*}
It is clear that the above process $v$ is valued in $\Gamma_1\times\Gamma_2$ since $\Gamma_1, \Gamma_2$ are cones. If $X_{t-}>0$, then
\begin{align*}
\phi(X,u)&=X^2\Big[v_1^{\top}(R_1+P_1D^{\top}D)v_1+2v_1^{\top}(P_1B_1+D^{\top}(P_1C+\Lambda_1))
-H_{11}^{\alpha_{t-},*}(P_1,\Lambda_1)\Big]\\
&\qquad+X^2\int_{\cZ}\Big[v_2^{\top}R_2v_2+2P_1B_2(z)^{\top}v_2(z)-2P_1\sum_{k=1}^{n_2}(E_k(z)X+F_k(z)v_2(z))\\ &\qquad\qquad+\sum_{k=1}^{n_2}(P_{1}+\Gamma_{1k})\Big(\big((1+E_k(z)+F_k(e)v_2(z))^+\big)^2-1\Big)\\
&\qquad\qquad+\sum_{k=1}^{n_2}(P_{2}+\Gamma_{2k})\big((1+E_k(z)+F_k(e)v_2(z))^-\big)^2\\
&\qquad\qquad-H_{12}^{\alpha_{t-},*}(P_1,P_2,\Gamma_1, \Gamma_2)\nu(\dz)\Big]\nu(\dz)\geq 0,
\end{align*}
from the definitions of $H_{11}^{i,*}, H_{12}^{i,*}$. Moreover, the equality holds at
\begin{align*}
u^*_1(t,X,\alpha)=\hat v_{11}^{\alpha_{t-}}(t,P^{\alpha_{t-}}_1,\Lambda^{\alpha_{t-}}_1)X_{t-}^+, \
u^*_2(t,X,\alpha)=\hat v_{12}^{\alpha_{t-}}(t,P^{\alpha_{t-}}_1,P^{\alpha_{t-}}_2,\Gamma^{\alpha_{t-}}_1, \Gamma^{\alpha_{t-}}_2)X_{t-}^+.
\end{align*}
Next if $X_{t-}<0$, then
\begin{align*}
\phi(X,u)&=X^2\Big[v_1^{\top}(R_1+P_1D^{\top}D)v_1-2v_1^{\top}(P_1B_1+D^{\top}(P_1C+\Lambda_1))
-H_{21}^{\alpha_{t-},*}(P_2,\Lambda_2)\Big]\\
&\qquad+X^2\int_{\cZ}\Big[v_2^{\top}R_2v_2-2P_2B_2(z)^{\top}v_2(z)-2P_2\sum_{k=1}^{n_2}(E_k(z)-F_k(z)v_2(z))\\ &\qquad\qquad+\sum_{k=1}^{n_2}(P_{1}+\Gamma_{1k})\big((-1-E_k(z)+F_k(e)v_2(z))^+\big)^2\\
&\qquad\qquad+\sum_{k=1}^{n_2}(P_{2}+\Gamma_{2k})\Big(\big((-1-E_k(z)+F_k(e)v_2(z))^-\big)^2-1\Big)\\
&\qquad\qquad-H_{22}^{\alpha_{t-},*}(P_1,P_2,\Gamma_1, \Gamma_2)\nu(\dz)\Big]\nu(\dz)\geq 0,
\end{align*}
from the definitions of $H_{21}^{i,*}, H_{22}^{i,*}$. Moreover, the equality holds at
\begin{align*}
u^*_1(t,X,\alpha)=\hat v_{21}^{\alpha_{t-}}(t,P^{\alpha_{t-}}_1,\Lambda^{\alpha_{t-}}_1)X_{t-}^-, \
u^*_2(t,X,\alpha)=\hat v_{22}^{\alpha_{t-}}(t,P^{\alpha_{t-}}_1,P^{\alpha_{t-}}_2,\Gamma^{\alpha_{t-}}_1, \Gamma^{\alpha_{t-}}_2)X_{t-}^-.
\end{align*}
Finally, when $X_{t-}=0$, then
\begin{align*}
\phi(X,u)&=u_1^{\top}R_1P_2D^{\top}Du_1+\int_{\cZ}\Big[u_2^{\top}R_2u_2
+\sum_{k=1}^{n_2}(P_{1}+\Gamma_{1k})\big((F_ku_2)^+\big)^2\\
&\qquad+\sum_{k=1}^{n_2}(P_{2}+\Gamma_{2k})\big((F_ku_2)^-\big)^2\Big]\nu(\dz)\geq0;
\end{align*}
here the equality holds at $u^*_1=0$, $u^*_2=0$.

The above analysis together with \eqref{lower2} shows that
\begin{align*}
&\E\Big[ P^{\alpha_{\tau_n}}_{1,{\tau_n}}(X_{\tau_n}^+)^2+P^{\alpha_{\tau_n}}_{2,{\tau_n}}(X_{\tau_n}^-)^2\Big]+\E\int_0^{\tau_n}\Big[u_{1}^{\top}R_{1}u_{1}
+QX^2+\int_{\cZ}u_{2}(z)^{\top}R_{2}(z)u_{2}(z)\nu(\dz)\Big]\dt\\
&\geq P^{i_0}_{1,0}(x^+)^2+P^{i_0}_{2,0}(x^-)^2.
\end{align*}

By noting that for any $u\in\mathcal{U}$, the corresponding state process $X\in S^{2}_{\mathbb{F}}(0,T;\R)$. Sending $n\rightarrow\infty$, we conclude, from the dominated convergence theorem, that \eqref{lower1} holds for any $u\in\mathcal{U}$, where the equality is achieved when $u^*$ is defined by \eqref{ustar}.
\eof
\begin{lemma}
Under the condition of Theorem \ref{Th:verif}, $(u^*_1(t,X,\alpha),u^*_2(t,X,\alpha))\in\mathcal{U}$.
\end{lemma}
\pf
It is clear that $(u^*_1(t,X,\alpha),u^*_2(t,X,\alpha))$ is valued in $\Pi_1\times\Pi_2$. It remains to prove $$(u^*_1(t,X,\alpha),u^*_2(t,X,\alpha))\in L^2_\mathbb{F}(0, T;\R^{m_1})\times L^{2,\nu}_{\mathcal{P}}(0, T;\R^{m_2}).$$

%
%
%

Substituting \eqref{ustar} into the state process \eqref{state}, we have
\begin{align}\label{stateoptimal}
\begin{cases}
\dd X_t=\left[AX+B^{\top}(\hat v_{11}X^++\hat v_{21}X^-)+\int_{\cZ}B(z)^{\top}(\hat v_{12}X^++\hat v_{22}X^-)\nu(\dz)\right]\dt\\
\qquad\qquad\qquad+\left[CX+D (\hat v_{11}X^++\hat v_{21}X^-)\right]^{\top}dw_t\\
\qquad\qquad\qquad+\int_{\cZ}\left[E(z)X+F(z)(\hat v_{12}X^++\hat v_{22}X^-)\right]^{\top}\tilde N(\dt,\dz), \ t\in[0,T], \\
X_0=x,\ \alpha_0=i_0,
\end{cases}
\end{align}
According to \cite[Theorem 3.5]{HSX}, we have $\hat v_{11}, \hat v_{21}\in L^{2}_{\mathbb{F}}(0, T;\R^{m_1})$. And from \eqref{boundeddomain}, we know that $\hat v_{12},\hat v_{22}\in L^{\infty}_{\mathcal{P}}(0, T;\R^{m_2})$. By the basic theorem of Gal'chuk \cite[p.756-757]{Gal}, the SDE \eqref{stateoptimal} admits a unique solution, denoted by $X^*$.

From the proof of Lemma \ref{lemma1}, we find that, for any stopping time $\iota\leq T$,
\begin{align}\label{iota}
&\E\Big[ P^{\alpha_{\theta_n}\wedge\iota}_{1,{\theta_n\wedge\iota}}((X^*_{\theta_n\wedge\iota})^+)^2
+P^{\alpha_{\theta_n}\wedge\iota}_{2,{\theta_n\wedge\iota}}((X^*_{\theta_n\wedge\iota})^-))^2\Big]\nn\\
&\qquad\qquad+\E\int_0^{\theta_n\wedge\iota}\Big[(u^*_{1})^{\top}R_{1}u^*_{1}
+Q(X^*)^2+\int_{\cZ}(u^*_{2}(z)^{\top}R_{2}(z)u^*_{2}(z))\nu(\dz)\Big]\dt\nn\\
&= P^{i_0}_{1,0}(x^+)^2+P^{i_0}_{2,0}(x^-)^2,
\end{align}
where
\begin{align*}
\theta_n:=\inf\{t\geq 0:|X^*_{t}|\geq n\}\wedge T.
\end{align*}

\underline{When Assumption \ref{assu2} holds.}
We have, from \eqref{iota},
\begin{align*}
\delta\E\int_0^{\theta_n\wedge T}\Big[|u^*_{1}|^2
+\int_{\cZ}|u^*_{2}(z)|^2\nu(\dz)\Big]\dt\leq P^{i_0}_{1,0}(x^+)^2+P^{i_0}_{2,0}(x^-)^2.
\end{align*}
Letting $n\rightarrow\infty$, it follows from the monotone convergence theorem that $(u^*_1,u^*_2)\in L^2_\mathbb{F}(0, T;\R^{m_1})\times L^{2,\nu}_{\mathcal{P}}(0, T;\R^{m_2})$.


\underline{When Assumption \ref{assu3} holds}. In this case, there exists $c>0$ such that $P^i_j\geq c, \ P^i_{j}+\Gamma^i_{jk}\geq c, \ i\in\cM, \ j=1,2, \ k=1,\ldots, n_2$. From \eqref{iota}, we get
\begin{align*}
c\E[ |X^*_{\theta_n\wedge\iota}|^2]\leq P^{i_0}_{1,0}(x^+)^2+P^{i_0}_{2,0}(x^-)^2.
\end{align*}
Letting $n\rightarrow\infty$, it follows from Fatou's lemma that
\begin{align*}
\E[ |X^*_{\iota}|^2]\leq c,
\end{align*}
for any stopping time $\iota\leq T$. This further implies
\begin{align*}
\E \int_{0}^{ T}|X^*_{t}|^2\dt \leq cT.
\end{align*}
Applying It\^o formula to $|X^*_{t}|^2$, yields that
\begin{align*}
&\quad x^2+\E\int_0^{\theta_n\wedge T} |Du^*_1|^2\dt+\E\int_0^{\theta_n\wedge T}\int_{\cZ}|Fu^*_2| ^2\nu(\dz)\dt\\
&= \E [X^*_{\theta_n\wedge T}]^2-\E\int_0^{\theta_n\wedge T}\Big[(2A+|C|^2+\int_{\cZ}|E(z)|^2\nu(\dz))|X^*_{t}|^2\\
&\qquad+2X^*_{t}(B_1+D^{\top}C)^{\top}u_1^* +2X^*_{t}\int_{\cZ}(B_2^{\top}u_2^*+\sum_{k=1}^{n_2}E_kF_ku_2^*)\nu(\dz) \Big]\dt.
\end{align*}

If  \ref{R12} and \ref{R22} hold, we have
\begin{align*}
&\quad \delta\E\int_0^{\theta_n\wedge T}|u^*_1|^2\dt+\delta\E\int_0^{\theta_n\wedge T}\int_{\cZ}|u^*_2| ^2\nu(\dz)\dt\\
&\leq c +c\E\int_0^{\theta_n\wedge T}\Big[|X^*_{t}|^2+2|X^*_{t}||u_1^* | +2|X^*_{t}|\int_{\cZ}|u_2^*|\nu(\dz) \Big]\dt\\
&\leq c +c\Big(1+\frac{2}{\delta}+\frac{2\nu(\cZ)}{\delta}\Big)\E\int_0^{\theta_n\wedge T}|X^*_{t}|^2\dt+\frac{\delta}{2}\E\int_0^{\theta_n\wedge T}|u_1^* |^2\dt +\frac{\delta}{2}\E\int_0^{\theta_n\wedge T}\int_{\cZ}|u_2^*|^2\nu(\dz) \dt.
\end{align*}
After rearrangement, it follows from the monotone convergence theorem that
\begin{align*}
\E\int_0^{ T}|u^*_1|^2\dt+\E\int_0^{ T}\int_{\cZ}|u^*_2| ^2\nu(\dz)\dt\leq c.
\end{align*}
Hence $(u^*_1,u^*_2)\in L^2_\mathbb{F}(0, T;\R^{m_1})\times L^{2,\nu}_{\mathcal{P}}(0, T;\R^{m_2})$.

If \ref{R12} and \ref{R21} hold, $u^*_1\in L^2_\mathbb{F}(0, T;\R^{m_1})$ follows exactly as above. On the other hand, we can get
$u_2^*\in L^{2,\nu}_{\mathcal{P}}(0, T;\R^{m_2})$ from \eqref{iota}.

The last case in Assumption \ref{assu3} can be handled similarly.
\eof

\end{document}